\magnification=\magstep 1
\input amstex
\documentstyle{amsppt}

\def\a{\alpha}
\def\b{\beta}
\def\ga{\gamma}
\def\O{\Omega}
\def\o{\omega}
\def\ti{\widetilde}
\def\D{\Delta}
\def\e{\epsilon}

\def\d{\delta}
\def\A{{\Cal A}}
\def\pa{\partial}
\def\np{{\hskip 0.8pt \hbox to 0pt{/\hss}\hskip -0.8pt\partial}}
\def\G{{\Cal G}}

\def\B{{\Cal B}}
\def\F{{\Cal F}}

\def\UU{{\Cal U}}

\def\DD{{\Cal D}}
\def\M{{\Cal M}}
\def\N{{\Cal N}}
\def\p{\phi}
\def\P{\Phi}
\def\R{{\Cal R}}
\def\si{\sigma}
\def\n{\nabla}
\def\t{\theta}
\def\lmd{\lambda}

\def\Ga{\Gamma}

\def\ZZ{\Bbb Z}
\def\X{{Cal X}}
\def\CC{\Bbb C}
\def\NN{\Bbb N}
\def\RR{\Bbb R}

\def\dt{\frac{d}{dt}}

\def\SS{{\Cal S}}
\def\W{{\Cal W}}
\def\SW{\hbox{\bf SW}}

\def\BA{[3]}
\def\Ba{[4]}
\def\BD{[6]}
\def\Bd{[7]}
\def\Fl{[8]}
\def\FU{[9]}
\def\Fu{[10]}
\def\Bg{[11]}
\def\KM{[12]}

\def\SA{[15]}
\def\Ta{[17]}

\def\WI{[19]}
\def\CCS{[5]}

\def\Bz{[16]}
\def\X{{\Cal X}}

\topmatter
\title   Bott-type and Equivariant Seiberg-Witten
Floer homology I  
\endtitle
(preliminary)
\author Guofang Wang and Rugang Ye
\endauthor
\address Max-Planck Institut f\"ur Mathematik in den Naturwissenschaften,
Inselstra{\ss}e 22-26, 04103 Leipzig, Germany
\endaddress
\email gwang\@mis.mpg.de
\endemail
\address Institute 
of System Science, Academia Sinica, 
100080 Beijing, China
\endaddress 
\email gfwang\@iss06.iss.ac.cn
\endemail
\address Department of Mathematics, University of
California, Santa Barbara, CA 93106, USA
\endaddress
\email yer\@math.ucsb.edu
\endemail
\address Ruhr-Universit\"at Bochum, Fakult\"at f\"ur Mathematik,
44780 Bochum, Germany
\endaddress
\email ye\@dgeo.ruhr-uni-bochum.de
\endemail

\rightheadtext{Bott-type  Seiberg-Witten
Floer homology I } 
\abstract
We construct Bott-type and stable equivariant 
Seiberg-Witten Floer homology 
and cohomology for rational homology spheres and prove their
diffeomorphism invariance.
\endabstract
\endtopmatter
\document

\head Table of Contents
\endhead
\noindent Part I 

\noindent 1. Introduction 

\noindent 2. Preliminaries

\noindent 3. Seiberg-Witten moduli space over $Y$

\noindent 4. Seiberg-Witten trajectories: transversality

\noindent 5. Index and orientation 

\noindent 6. The temporal model and compactification

\noindent 7. Bott-type 
homology

\noindent 8. Invariance I

\noindent Appendix A  Seiberg-Witten 
Floer homology

\noindent Appendix B. Equivariant Seiberg-Witten 
Floer homology

\noindent Appendix C. Two analysis lemmas 
\bigskip
\noindent From Part II

\noindent 9. Stable Bott-type and  stable equivariant 
homology  

\noindent 10. Invariance II





\head 1. Introduction
\endhead

At the very beginning of the development of the new 
Seiberg-Witten gauge theory it was clear that,
at least formally, the celebrated instanton homology theory of
A. Floer for 3-manifolds (homology spheres) $\Fl$ could be adapted to the
Seiberg-Witten set-up. Indeed, the original
4-dimensional Seiberg-Witten equation leads naturally
to a 3-dimensional Seiberg-Witten equation via a
limit process, as first observed by Kronheimer and
Mrowka $\KM$. To establish a Seiberg-Witten Floer
homology for a 3-manifold $Y$, the obvious idea is to replace flat 
connections in Floer's set-up by solutions of the 3-dimensional 
Seiberg-Witten equation on $Y$ (henceforth called {\it Seiberg-Witten
points}), and instanton trajectories by Seiberg-Witten
trajectories, which are solutions of the 4-dimensional Seiberg-Witten 
equation on the infinite cylinder $Y \times \RR$. Note that the 
Seiberg-Witten points are precisely the critical 
points of the Seiberg-Witten type Chern-Simons functional, and 
that the Seiberg-Witten trajectories are 
precisely the trajectories (negative gradient flow 
lines) of this functional.  Hence the said idea amounts to establishing 
a Morse-Floer theory for the Seiberg-Witten type Chern-Simons functional.
However, one encounters various difficulties
when trying to implement this idea. One most
serious problem is that Seiberg-Witten Floer homologies for a homology
sphere (or rational homology sphere) may depend on the underlying 
Riemannian
metric (cf. e.g. $\Bd$), and hence are generally not diffeomorphism
invariants. The purpose of this paper is to resolve this problem.

The trouble, e.g. in the situation of homology spheres, 
is caused by the trivial Seiberg-Witten point: the trivial connection
coupled with the zero spinor field. It is reducible,
i.e. it is a fixed point of the action of the group $S^1$ of constant 
gauges.
Under reasonable perturbations of the Seiberg-Witten equation, this
reducible point always survives. To deal with it, one can use suitable
perturbations to make it a transversal point for
the Seiberg-Witten equation. Then one can
construct a Seiberg-Witten Floer homology, see Appendix A.
However, one encounters a serious obstruction when
trying to compare the homologies for two different perturbation parameters
(e.g. metrics). A canonical way of such comparison is to construct 
chain maps in terms of parameter-dependent Seiberg-Witten trajectories 
which connect the 3-dimensional Seiberg-Witten equation of one parameter 
to that of another. We shall call them {\it transition 
trajectories.} The said obstruction is the  presence of reducible 
transition trajectories with negative spectral flow of the
linearized Seiberg-Witten operator. Such trajectories are not in 
transversal 
position and may appear in the compactification 
of the moduli spaces of transition trajectories between 
irreducible Seiberg-Witten points. Consequently, 
the compactified moduli spaces of transition trajectories may be very 
pathological and cannot be used to define the desired chain maps. 

The appearance of such trajectories roughly goes as follows.
The spectral flow along a reducible Seiberg-Witten
trajectory for a fixed parameter is $1$. That along a reducible 
transition trajectory from a given generic 
parameter to a nearby one is also 1. When passing from one 
generic parameter to another through certain degenerate parameters,
the spectral flow jumps and becomes negative.  
Here, typically,  Seiberg-Witten Floer homology also jumps. 

\bigskip 
\noindent {\bf Bott-type Construction}
\smallskip
Since the ``ordinary" Seiberg-Witten Floer homologies
may not be diffeomorphism invariants for rational
homology spheres, we seek an alternative construction.
It turns out that the Bott-type set-up above the
$S^1$ quotient is a right one. In other words, we work on the level of
quotient by the group of based gauges rather
than the full group of gauges. {\it But even in this set-up,
one has to choose the right approach in order to 
obtain an invariant theory.} 

Now, in this alternative set-up, the irreducible part of the moduli space 
of gauge classes of Seiberg-Witten points consists of 
finitely many circles, while its reducible part consists of a single point, 
provided that we choose a generic parameter. These circles and the 
reducible 
point are precisely the critical submanifolds of the (Seiberg-Witten
type) Chern-Simons functional. Our goal here amounts to establishing 
a Bott-type Morse-Floer theory for the Chern-Simons
functional. The basic strategy is to use the moduli spaces of 
trajectories between critical submanifolds to send
(co)homological chains from one critical submanifold to others, which  
defines the desired boundary operator
for the (co)chain complex. This is a natural
extension of Floer's construction and was first used by Austin-Braam 
$\BA$ and Fukaya $\FU$ in Floer's set-up. The former authors use 
differential forms as chains and cochains, while the 
latter uses ``geometric chains". We shall adopt the classical singular
chains and cochains.  We emphasize that {\it it is not 
clear whether the differential form approach leads  to  diffeomorphism
invariants.}  On the other hand, an advantage of the singular chain 
set-up is that it admits integer coefficients and hence 
contains torsion information.

\bigskip
\noindent {\bf Stable Equivariant Construction}
\smallskip

On the level of the based gauge quotient, the group 
$S^1$ of constant gauges acts and equivariant 
Seiberg-Witten Floer homologies can be defined. For
example, one can follow the approach of $\BA$, 
see Appendix B.  We 
can also use equivariant singular cochains, see
Appendix B.  But we do not yet know whether these 
homologies and cohomologies are diffeomorphism 
invariants.  (As explained before, the essential 
trouble is caused by reducible transition  trajectories with negative 
spectral flow.  If they are present, {\it 
in general the traditional construction of the desired chain maps 
breaks down completely.} Of course, this does not 
exclude the possibility that there might be other undiscovered 
schemes for constructing  chain maps.
{\it But we tend to believe that equivariant 
Seiberg-Witten Floer homologies are not diffeomorphism invariants.}) 
 Instead we construct a ``stable" equivariant cohomology which we 
show to be invariant. The word ``stable" means that 
we couple the configuration space 
of the Seiberg-Witten  equation with the unit circle,
hence increasing the dimension of the moduli 
space of Seiberg-Witten points by one. In this 
set-up, we first obtain a stable Bott-type theory
in the same way as the above Bott-type theory.
Restricting  to equivariant singular cochains, we 
then obtain the stable equivariant Seiberg-Witten 
Floer homology and cohomology.  Note that one can also 
define a stable equivariant homology using differential forms analogous 
to the differential form construction in Appendix B, 
but  it is not clear whether it is 
an invariant.  Hence both in the Bott type 
and stable equivariant cases, {\it it is important 
to use the approach of singular chains (or similar objects) instead of 
differential forms.} 

\bigskip
\noindent {\bf Spinor Perturbation}
\smallskip
To prove the diffeomorphism invariance of the Bott-type
Seiberg-Witten Floer (co)homology, we have
to overcome the said obstruction of reducible
transition trajectories with negative spectral flow. Our first strategy 
for this is to perturb the spinor equation in
the transition trajectory equation 
in order to eliminate these transition 
trajectories. We utilize the vanishing of
the rational homology group to construct suitable
vector fields which are equivariant under
based gauges. Note that {\it they are not equivariant 
under constant gauges}.  A desired perturbation is then gotten by adding 
one of these vector fields to 
the spinor equation. Thus the source of our trouble,
namely the vanishing of the rational homology
group, also works to our benefit - a rather amusing phenomenon. 
There is another aspect here: in one step of the 
invariance proof,
we have to show that the reducible
trajectory from the reducible Seiberg-Witten point
to itself for a fixed parameter is transversal. Here, we again use 
the vanishing of the rational homology group. (One can also use the 
Fredholm perturbations described below, which are however more 
complicated.) 

The spinor perturbations will be applied 
to the stable set-up as well. 
Here {\it they are also equivariant under constant 
gauges}. With these perturbations, we can immediately 
prove the diffeomorphism invariance of the stable Bott-type theory and
the stable equivariant theory. 

However, for the Bott-type set-up without stablization, 
a second strategy is also needed. Indeed, here the spinor perturbation will 
be applied whenever the spectral flow along reducible
transition trajectories 
from one given parameter $t_1$ to another
$t_2$ is nonpositive. If the spectral flow 
equals $1$, we perform the spinor perturbation
only in one direction between $t_1$ and $t_2$,
say from $t_1$ to $t_2$. In either case, the spectral flow 
in the reversed direction from $t_2$ to $t_1$ is at least $1$.  
We need to achieve transversality in this direction as well, but 
here we are no longer allowed to use the spinor
perturbation. The reason is that if we 
apply the spinor perturbation in both directions, the gauge invariance 
group for the glued equation from $t_1$ to $t_1$ would be too small. 

\bigskip  
\noindent {\bf Cokernel Perturbation}
\smallskip
Our second strategy is to use Kuranishi models near 
reducible transition trajectories to achieve 
transversal perturbations which are equivariant under the full gauge 
group. These perturbations are given in terms of suitable operators onto  
the cokernel of the linearized Seiberg-Witten operator, 
hence we call them {\it cokernel perturbations}.   Such perturbations 
were first used by Donaldson $\BD$  in a geometric context. 
Note that the domain of our cokernel perturbations 
is the kernel of the linearized Seiberg-Witten operator. 
This feature depends  on the positive spectral flow condition. 
(If the spectral flow is negative, cokernel perturbations can also be 
performed to achieve transversality, but 
the equivariance can only be preserved, provided that 
additional parameters are introduced.
By adding additional parameters, however, we arrive at the stable set-up. 
The resulting perturbations are more complicated than 
the spinor perturbations, hence we chose to use the latter for 
the stable set-up.)

After performing the cokernel perturbations,
the moduli space of transition trajectories between the reducibles 
becomes transversal. To achieve transversality for the 
compactified moduli spaces of transition trajectories 
between irreducibles or one reducible and one 
irreducible, we have to find a way to extend these 
perturbations to the situation of gluing. In other words, we have to 
perform cokernel perturbations near infinity, i.e.  
the boundary of the compactified moduli spaces.
We employ extended Kuranishi models, or Kuranishi models around glued 
approximate transition trajectories to construct 
the desired perturbations.  Then we can glue reducible transition 
trajectories with trajectories and obtain the desired transversality 
for the compactified moduli spaces. This 
important construction will be used in several steps.  

\bigskip
\noindent{\bf Further Analysis}
\smallskip

There are a few further delicate analytical 
issues in the above constructions we would like 
to address here.  First, as explained before, we use moduli spaces 
of (transition) trajectories to define our boundary 
operator and establish the equivalence isomorphism. 
There is a subtlety here in the choice of the moduli spaces.
Namely we need appropriate endpoint maps from the moduli 
spaces of trajectories to the moduli spaces of Seiberg-Witten 
points.  For this purpose,  we use the temporal model for the trajectory 
spaces. On the other hand, we have to show that the compactified 
moduli spaces have a structure of smooth manifolds with 
corners. It is not clear how to prove this directly for the temporal model.
We work instead with another model, and introduce a "twisted 
time translation" action  on it, which is induced from the time 
translation action on the temporal model. Even for this model, 
considerable care is needed for establishing the structure of smooth 
manifolds with corners.

A fundamental issue here is convergence (modulo
splitting) of trajectories. When splitting
occurs, we have to show the important property that the
intermediate endpoints of the limit trajectories match each
other. This requires a strong and detailed convergence
result (Proposition 6.8).

The spinor perturbations cause additional analytical difficulties which 
demand special treatments. For example, one has to establish 
a uniform $L^{\infty}$ estimate for the spinor part of the 
parameter-dependent Seiberg-Witten trajectories. 
With the presence of the perturbations, the ordinary pointwise maximum 
principle argument does not work. Instead, we apply the 3-dimensional 
Weitzenb{\"o}ck formula (rather than the 4-dimensional one)
to obtain an initial local integral estimate in terms of 
the Seiberg-Witten energy. Then we  apply the 4-dimensional 
Weitzenb{\"o}ck formula and the technique of Moser iteration to 
derive the desired $L^{\infty}$ estimate.

\bigskip
Now we make a few concluding remarks. First, using the invariants 
constructed in this paper and Seiberg-Witten 
Floer homology for manifolds with nonzero first Betti number, 
one can define  relative Seiberg-Witten invariants 
for general four dimensional manifolds with boundary.
This will be discussed elsewhere.
(Manifolds with first betti number equal to one requires special 
treatment.) Second, we would like to mention that the theory in 
this paper can be strengthened to a Seiberg-Witten Floer homotopy  theory 
along the lines of the Floer homotopy theory as proposed in $\CCS$. 
which amounts to the ordinary 
Seiberg-Witten theory  with transversality replaced by Bott-type 
transversality. If we lift to the level of the based gauge 
quotient, we obtain the double Bott-type 
Seiberg-Witten Floer homology, which is isomorphic
to the Bott-type Seiberg-Witten Floer homology.
Similarly, if the manifold $Y$ has a symmetry group, then we can consider 
the equivariant theory with respect to the group. Lifting 
it to the stable set-up, we obtain a double 
equivariant theory. Details will appear elsewhere.

\bigskip
A major part of the results in this work 
were obtained in Spring  1996 while both authors were at Bochum 
University. 
This work has been reported by the second
named author in a number of talks given in 1996.

This work consists of two parts. The present 
paper is a preliminary version of Part I 
in which Sections 9 and 10 of Part II are also included. 


\head 2. Preliminaries
\endhead

To fix notations, we first recall the definitions
of the Seiberg-Witten equations on 3 and 4 dimensional manifolds. 

Let $(X_0, g)$ be an oriented Riemannian manifold of dimension
$n$ and $Spin^c(X_0)$ the set of isomorphism classses
of $spin^c$ structures on $X_0$.
Consider a $spin^c$ structure $c \in Spin^c(X_0) $ and 
its associated spinor bundle $W$ and line bundle $L$. (More
precisely, $c$ is a representative of an element in $Spin^c(X_0)$.
The homology invariants we are going to construct depend 
are {BBBB depend are ??BBBB} independent of the choice of the 
representative.)
We have the associated configuration space $\A \times \Ga(W)$,
where $\A$ denotes the space of smooth unitary connections 
on $L$ and $\Ga$ the 
space of smooth sections of a vector bundle. (We suppress 
the dependence on $c $ in the notations.) The 
gauge transformation group (the group of 
gauges) is $\G = C^{\infty}(X_0, S^1)$, where $S^1 \equiv U(1)$ denotes
the unit circle in $\CC$. ($\G$ depends only on $X_0$.)
The action of $\G$ on $\A \times \Ga(W)$
is defined by $((A, \P), g) \to (A + g^{-1}dg, g^{-1}\P)$.  
This formula also defines 
the (separate) actions of $\G$ on $\A$ and $\Ga(W)$. $\G$
acts freely on the subspace
of pairs $(A, \P)$ with $\P \not \equiv 0$. Such pairs
are called irreducible. The isotropy subgroup  at any reducible
pair $(A, 0)$ is the subgroup of gauge transformations 
which are constants on each component of $X_0$. If $X_0$ is 
connected, we identify it with 
$S^1$. In this case, we fix a reference point
$x_0 \in X_0$ and set $\G^0 = \{g \in \G : g(x_0)= 1\}$, which is
called the group of based gauges. Then
the quotient $\G \slash S^1$ is represented by $\G^0$. 

The action of a gauge $g$ will be denoted by 
$g^*$. We set $\B = (\A \times \Ga)  \slash \G$ and $\B^* = (\A \times
(\Ga -\{0\}))\slash \G$.  Let $\O^k(X_0)$ denote 
the space of smooth imaginary valued $k$-forms, and 
$\O^+(X_0)$ the space of smooth imaginary 
valued self-dual 2-forms (in the case that $\hbox{ dim }X_0 =4$).

We shall need the following 
\proclaim{Lemma 2.1}  
Assume that $X_0$ is closed.  Then the 
map from $\G$ to $H^1(X_0;\ZZ) \slash \{torsions\}$ given by $g \to 
$ the deRham class of $g^{-1}dg$ is 
surjective and induces an isomorphism 
from the component group of $\G$ to $H^1(X_0;\ZZ) \slash \{torsions\}$.
Moreover, there is a ~unique harmonic map $g$ with $g(x_0) = 1$ in each
component of $\G$, provided that 
$X_0$ is connected and $x_0 \in X_0$ is a fixed point.  
In particular, $\G$ is connected if $X_0$ is connected and
$H^1(X_0;\ZZ)$ is torsion.
\endproclaim

\demo {Proof} For simplicity,
assume that 
$X_0$ is connected.  The surjectivity of the said map
follows from integration along paths.  If $g^{-1}dg $
and $g_1^{-1}dg_1$ represent the same cohomology class,
then $g_1 = g e^{f}$ for some $f \in \O^0(X_0)$
as one easily sees. Hence $g_1$ and $g$ lie in the 
same component group. 

The statement about harmonic representative follows
from the standard theory of harmonic maps. It can 
also be derived quickly in an elementary way. For
example, if $g_1 = g e^f$ and $g$ are two 
harmonic maps, then $f$ is a harmonic function, hence
constant. \qed
\enddemo

We continue with 
the above $spin^c$ structure $c$ on $X_0$. A connection $A \in \A$ induces 
along with the Levi-Civita connection 
a connection $\n^A$ on the spinor bundle $W$ and 
the associated Dirac operator 
$D_A: \Gamma (W)
\to\Ga(W)$, 
$$D_A=\sum^n_{i=1}e_i\cdot\n^A_{e_i},$$
where $\{e_i\}$ denotes a
local orthonormal tangent frame and the 
dot denotes the Clifford multiplication. The Dirac operator is gauge 
equivariant, i.e. $D_A (g^{-1} \P) = g^{-1}D_A \P$,
and satisfies 
the following fundamental Weitzenb\"ock formula for the
Dirac operator
$$D^*_AD_A\P=-\Delta _A\P+\frac{s}{4}\P-\frac{1}{2}F_A\cdot\P,
\tag 2.1
$$
where $s$ denotes the scalar curvature of $(X,g)$ and   
$F_A$ the curvature of $A$.

Now we specify to the dimension $n = 4$. 
There is a canonical decomposition $W = W^+ \oplus W^-$ of the spinor 
bundle $W$.
The Dirac operator ~splits: $D_A : \Ga (W^+) \to \Ga (W^-)$,
$D_A: \Ga (W^-) \to \Ga (W^+)$. For a positive spinor field $\P \in
\Ga (W^+)$, the curvature $F_A$ in the above Weitzenb{\"o}ck formula
reduces to its self-dual part  $F_A^+$.

\definition{Definition 2.2}
The Seiberg-Witten equation with 
the given $spin^c$ structure $c$  is
$$
\eqalign{
F^+_A &=\frac{1}{4}\langle e_ie_j\P,\P\rangle e^i\wedge e^j,\cr
D_A\P &=0, \cr}
\tag 2.2
$$
for $(A, \P) \in \A \times \Ga(W^+)$, where 
$\{e^i\}$ denotes the
dual of $\{e_i\} $ (a local orthonormal tangent frame).
The Seiberg-Witten operator is
$$
\hbox{\bf SW}(A, \P) = (F^+_A -\frac{1}{4}\langle e_ie_j\P,
\P\rangle e^i\wedge e^j,  D_A\P).
$$
\enddefinition

Note that the Seiberg-Witten 
operator is gauge equivariant, where $\G$ acts on 2-forms trivially. 
Consequently, the Seiberg-Witten equation is gauge invariant.  

Next let $(Y,h)$ be an oriented, closed Riemannian 3-manifold with metric
$h$, and $c$ a $spin^c$ structure on $Y$. We have the associated spinor 
bundle  $S= S_c(Y)$, line bundle $L(Y) = L_c(Y)$ and   
the other associated spaces: $\G(Y)$, $\A(Y)$, $\B(Y)$, $\B^*(Y)$ etc.
We set $X = Y \times \RR$, which will be equipped with the product 
metric and given the  orientation $(e_1, e_2, e_3, \frac{d}{dt})$, where
$(e_1, e_2, e_3)$ denotes a positive local orthonormal frame on $Y$. Let 
$\pi : X \to Y$ denote the projection. 
The $spin^c$ structure $c$ induces a $spin^c$ 
structure $\pi^* c$ on $X$  with the associated line bundle $L_X =
\pi^* L_Y$ and associated spinor bundles $W^+ = \pi^* S$ 
, $W^-$.  We have the following relation between the 
Clifford multiplications on $S$ and on $W^+$:
$$v\cdot \p (y) = -(\dt\cdot v\cdot \pi^* \p) (y, 0).$$

The associated spaces for $X$ will be indicated by 
the letter $X$. 
Let $i_t: Y \to X$ denote the inclusion map which sends $y \in Y$ to
$(y, t) \in X$. A connection $A \in \A(X)$  can be written as 
$$A=a(t)+f(\cdot, t)dt,$$
where $a(t)=i_t^*(A) \in \A(Y)$ and $f \in C^{\infty}(X, i\RR)$.
We set $\p(t)= \P(\cdot, t)=
i^*_t(\P)$ for $\P \in \Ga(W^+)$.
With these notations, we can rewrite (2.1) as follows
$$
\eqalign{
\frac{\partial a}{\partial t} &=*F_a+d_Yf+ \langle e_i\cdot\p,\p\rangle 
e^i,
\cr
\frac{\partial \p}{\partial t} &=-\np_a\p-f\p.\cr}
\tag 2.3
$$
Here $F_a$ denotes the curvature of $a$, $*$ the Hodge star operator
w.r.t. $h$, $d_Y$ the exterior
differential on $Y$, and $\np_a$ the Dirac operator associated with 
the connection $a$. 

\definition{Definition 2.3}
The Seiberg-Witten energy of $(A, \P)$ 
is
$$\eqalign{
E(A, \P) 
& = \int_X(|\frac{\pa\p}{\pa t}+f\p|^2+|\np \p|^2+|\frac{\pa a}
{\pa t}-df|^2+
|*F_a+\langle e_i\cdot\p,\p\rangle e^i|^2).\cr
}
$$
(The volume form is omitted.)
One readily shows that it is gauge invariant.
\enddefinition 

Using the finite energy condition one easily derives from  (2.3) 
the following limiting equation for a connection 
$a\in \A(Y)$ and a spinor field $\p \in \Ga(S)$
$$
\eqalign{
*F_a+\langle e_i\cdot\p,\p\rangle e^i& =0,\cr
\np_a\p & =0.\cr}
\tag 2.4
$$

\definition{Definition 2.4} The Seiberg-Witten equation 
on $Y$ with the $spin^c$ structure $c$ is defined to be (2.4).
The Seiberg-Witten operator on $Y$ is
$$
\hbox{\bf sw}(a, \p) = (*F_a+\langle e_i\cdot\p, 
\p\rangle e^i,
-\np_a\p). 
$$ 

\enddefinition

As in dimension 4, the Seiberg-Witten operator $\hbox{\bf sw}$ is 
gauge equivariant ($\G_Y$ acts trivially on 1-forms). 

We shall need the following perturbed 
Seiberg-Witten equation

$$
\eqalign{
*F_a+\langle e_i\cdot\p,\p\rangle e^i& = \nabla H(a),\cr
\np_a\p + \lambda \p & =0,\cr}
\tag 2.5
$$
where  
$\lambda$ denotes a real number  and 
$\nabla H$ the  $L^2$-gradient of a $\G(Y)$-invariant 
real valued function $H $ on $\A(Y) $. (The 
$L^2$-product is given in (2.7) below.)
We have the associated perturbed 
Seiberg-Witten operator $\hbox{\bf sw}_{\lambda, H}$.
Note that $\nabla H$ is gauge equivariant and belongs
to ${ker~ } d^*$, which are consequences of 
the gauge invariance of $H$.

The classical Chern-Simons functional plays 
a fundamental role in Floer's instanton homology 
theory. Similarly, a Chern-Simons functional associated 
with the 3-dimensional Seiberg-Witten equation will 
be important in our situation. This
functional was first used by Kronheimer and Mrowka in  their
proof of the Thom conjecture $\KM$.

\definition{Definition 2.5} The Chern-Simons functional with
respect to a reference connection $a_0$ is
$$\hbox {\bf cs}(a,\p)=\frac12\int_Y(a-a_0)\wedge (F_a+F_{a_0})+
\int_Y\langle\p,\np_a\p
\rangle.$$

Let  $\lambda$ and $H$ be as above.
The perturbed Chern-Simons functional with perturbation
$(\lambda, H)$ is

$$
\eqalign{
\hbox {\bf cs}_{\lambda,  H}(a,\p)=\frac12\int_Y(a-a_0)\wedge (F_a+
F_{a_0} )+
\int_Y\langle\p,\np_a\p\rangle  \cr
 - \lmd\int_Y\langle\p,\p\rangle 
+H(a, \p). \cr}
$$

\enddefinition

Under a gauge $g$
 the perturbed Chern-Simons functional
changes as follows:
$$\hbox {\bf cs}_{\lambda, H}(g^*(a,\p))=\hbox{\bf cs}_{(
\lambda, H)}(a,\p)+
2\pi i\int_Y c_1(L(Y))\wedge g^{-1}dg.
\tag 2.6$$
This formula implies that ${cs}_{\lambda, H}$
is invariant under the identity component of
$\G(Y)$. Hence it descends to  the quotient $\B(Y)$, 
provided that $Y$ is a rational homology sphere.

We introduce an $L^2$-product 
on $\O^1(Y) \oplus \Ga(S)$:

$$
\langle(\p_1,a_1),(\p_2,a_2)\rangle_{L^2}=\int
_Y(Re
\langle\p_1,\p_2\rangle+\langle a_1,a_2\rangle)
\tag2.7
$$
Here, $\langle a_1,a_2\rangle$ denotes the (pointwise) Hermitian product.
Easy computations lead to

\proclaim{Lemma 2.6}  The $L^2$-gradient of the perturbed 
Chern-Simons functional is given by
$$ \nabla \hbox {\bf cs}_{\lambda, H} (a, \p) = -\hbox{\bf sw}_{\lambda, H}.$$
It follows that the critical points of the 
perturbed Chern-Simons functional are precisely 
the solutions of the perturbed Seiberg-Witten 
equation.
\endproclaim

Consider a solution $(A, \P)$ of the Seiberg-Witten equation 
on the product $X$. Using a suitable 
 gauge we can transform it into temporal
form. Let's assume that it is already in
temporal form, i.e. 
$f \equiv 0$ in the formula $A = a + fdt$. 
Then Lemma 2.6 and the equation (2.3) imply
($\p(t)= \P(\cdot, t)$) 
$$
\frac{\partial }{\partial t} (a, \p) = - \nabla \hbox
{\bf cs}_{\lambda, H}(a, \p).
\tag 2.8
$$

Hence 
solutions of the Seiberg-Witten equation on the product 
$X$ can be interpreted as trajectories (negative gradient 
flow lines) of the Chern-Simons functional. A similar 
formula and statement hold for solutions  of 
the perturbed Seiberg-Witten equation on $X$, 
which is  
$$
\eqalign{
F^+_A =&\frac{1}{4}\langle e_ie_j\P,\P\rangle e^i\wedge e^j
\cr 
&+  \nabla H(a)\wedge dt+*(  \nabla H(a) \wedge dt),\cr
D_A\P =& - \lambda \frac{d}{dt} \cdot \P,\cr}
\tag 2.9
$$ 
where $\frac{d}{dt} \cdot
\P$ denotes the Clifford multiplication on $X$, and $A= a + fdt 
, \p(t) = \P(\cdot, t)$ as before. The operator 
$\hbox{\bf SW}_{\lambda, H}$ is defined in an obvious way.

Next we introduce the perturbed Seiberg-Witten 
energy : 

$$
\eqalign{
E_{\lambda, H}(A, \P) = 
\int_X(|\frac{\pa\p}{\pa t}+f\p|^2+|\np \p + \lambda \p|^2+|\frac{\pa a}
{\pa t}-df|^2 \cr
+|*F_a+\langle e_i\cdot\p,\p\rangle e^i - \nabla H(a, \p)|^2). \cr}
\tag 2.10
$$

Note that it is invariant under the  action of $\G(X)$.

\proclaim{Lemma 2.7} Assume that $ Y$ is a rational homology sphere. 
Let $A=a+fdt$  and the gauge 
equivalence class of $(a,\p)$ converges to $\a,\b\in \B(Y)$ as
$t\to-\infty,\infty$ respectively, then we have
$$\eqalign{E_{\lambda, H}(A, \P) =&
2\hbox{\bf cs}_{\lambda, H}(\a)
-2\hbox{\bf cs}_{\lambda, H}(\b)+
\int_X|D_A\P|^2\cr
&+2\int_X|
F^+_A-\frac{1}{4}\langle e_ie_j\P,\P\rangle e^i\wedge e^j
- dt \wedge \nabla H(a) -*(dt \wedge  \nabla H(a))|^2.\cr}\tag 2.10
$$
In particular, there holds for a solution 
$(A, \P)$ of (2.8) 

$$E_{\lambda, H} (A, \P)=2\hbox{\bf cs}_{\lambda, H}(\a)
-2\hbox{\bf cs}_{\lambda, H}(\b).
$$
\endproclaim

\demo{Proof}There is a similar computation in $\SA$.
We have
$$\eqalign{  \int_X|D_A\P|^2
+2\int_X|
F^+_A-\frac{1}{4}\langle e_ie_j\P,\P\rangle e^i\wedge e^j
- dt \wedge \nabla H(a) -*(dt \wedge  \nabla H(a))|^2 \cr
=\int_X(|\frac{\pa\p}{\pa t}+f\p+\np \p + \lambda \p|^2+|-\frac{\pa a}
{\pa t}+df+ 
*F_a+\langle e_i\cdot\p,\p\rangle e^i - \nabla H(a, \p)|^2) \cr
=E(\P,A)+2\!\!\int\!\langle \frac{\pa\p}{\pa t}\!+f\p,\np_a\p+\lmd\p\rangle
+\!2\langle -\frac{\pa a}
{\pa t}+df,*F_a \! 
+\!\langle e_i\cdot\p,\p\rangle e^i\! - \!\nabla H(a, \p)\rangle\cr
=2\dt\int(|\p|^2+\langle a-a_o,F_a\rangle)+2H(a)+2
\int(\langle \frac{\pa\p}{\pa t},\np_a\p\rangle+\langle\frac{\pa a}{\pa t}\p,
\p\rangle ).\cr}
$$
Since the metric on $X=Y\times\RR$ is the product metric,
we have
$$\dt\np_a\p-\np_a\dt\p= \frac{\pa a}{\pa t}\p. \tag 2.11$$
The desired conclusion follows.\qed\enddemo

\head 3. Seiberg-Witten moduli spaces over $Y$
\endhead

We continue with the $(Y, h)$ and $c$ of the last section. 
While our theory applies 
to arbitrary closed $Y$, we 
assume for convenience that $Y$ is connected. 
Fix a reference connection $a_0$. If $L_Y$ is a trivial bundle, we 
choose $a_0$ to be the
trivial  connection.  
We have $\A(Y)=a_0+\O^1(Y)$.
We shall use the $(l,p)$-Sobolev norms (the $L^{l, p}$-norms) 
for $l \geq 0$ and $p >0$:
$$\|u\|_{l,p}=(\sum_{0\le k\le l}\int_Y|\n^ku|^p)^{1/p}.$$
Consider the Sobolev spaces $\A_{l,p}(Y)$ and  $\Ga_{l,p}(S)$, which 
are the completions of $\A(Y)$ and $\Ga(S)$ with respect to  the 
$(l,p)$-Sobolev norm respectively. Similarly, we have  
the Sobolev spaces $\Omega^k_{l,p}(Y)$. The corresponding group of gauges is 
$\G_{l+1,p}(Y)$,
which is the completion of $\G(Y)$ with respect to the 
$(l+1, p)$-Sobolev norm.

We need to make a choice of the configuration spaces
$\A_{l, p}(Y) \times \Ga_{l, p}(S)$. We require $3p/(3-lp)
>3$ or $lp > 3$, for then all elements in $\A_{l, p}(Y) 
\times \Ga_{l,p}(S)$ are continuous, and hence the 
holonomy perturbations in the sequel can be performed. Moreover,
the corresponding gauges on the product space $X = Y
\times \RR$ are continuous. In particular, if we choose 
$l=1$, then we require $p >3$.

\definition {Definition 3.1} We have the following spaces of
Seiberg-Witten points
$$\SS\W_{l,p} = \SS \W_{h,\lambda, H, l, p} =\{(a, \p)
\in \A_{l, p}(Y) \times \Ga_{l, p}(S): \hbox{\bf sw}_{\lambda, H}((a, \p))
=0 \}
$$
and the following moduli spaces of Seiberg-Witten points
$$\eqalign{\R_{l, p}  &=\R_{h, \lambda, H, l, p} = {\Cal SW}_{l,p}\slash
\G_{l+1, p}(Y), \cr
\R^0_{l,p} &=\R^0_{h, \lambda, H, l, p} = {\Cal SW}_{l,p} \slash \G^0_{l+1,
p}.\cr}
$$
The irreducible part of e.g. $\R_{l, p}$ will be denoted by $\R^*_{l,
p}.$
\enddefinition

We choose to work with the configuration space $\A_{1, 4}(Y)
\times \Ga_{1 ,4}(Y)$. Henceforth, the subscript $l$ stands for 
$(l, 4)$, e.g. $\|u\|_l = \|u\|_{l, 4},\A_l = \A_{l, 4}.$ 
We have the quotients $\B_{1}(Y)$ and 
$\B^0_{1}(Y)$ of the chosen configuration 
space $\A_1(Y) \times \Ga_1(Y)$ under the group of 
gauges $\G_2(Y)$ and 
group of based gauges $\G_2^0$. The gauge class of 
$(a, \p)\in \A_{1}(Y) \times \Ga_{1}(S)$ with 
respect to the full gauge group $\G_2(Y)$ will be 
denoted by $[a, \p]$. 
Its gauge class with respect to based gauges
will be denoted by $[a, \p]_0$.

\proclaim{Lemma 3.2} 1)Each element in $\R_{1}$ or
$\R^0_{1}$ can be represented 
by a smooth pair $(a, \p)$.

2) If $\|\n H\|_{L^{\infty}} < C$ for a constant $C$, then $\R_{1}$ 
and $\R^0_{1}$ are compact.
\endproclaim

\demo{Proof} We present the  proof for 2), which contains 
the argument for 1). Let $(a, \p) \in \A_1(Y) \times \Ga_1(S)$ be 
a solution of (2.5). Applying the 3-dimensional  Weitzenb\"ock formula in 
the weak form, the bound  $\|\n H\|_{
L^{\infty}}<C$ and Moser's weak maximum principle (cf. the 
proof of Proposition 8.5),  we obtain 
 $\|\p\|_{L^{\infty}}<C$ for a constant $C $. (Here and in the sequel, 
 we use the same letter $C$ to denote all constants which appear 
in a priori estimation).
Since $H^1(Y,\RR)=0$, by Hodge decomposition, $a-a_0=d\ga+d^*\d$ for
some $\ga\in\O_{0}^0(Y)$ and $\d\in\O_{0}^2(Y)$.  
 By gauge fixing, we can assume $d^*(a-a_0)=0$. Note that 
we can achieve this gauge fixing by 
a based gauge. Hence $a-a_0=d^*\d$. Furthermore, we can assume
that $d\d=0$. Hence we have $\D\d=F_a-F_{a_0}$. Since $\|F_a-F_{a_0}\|\le
C\|\p\|^2+C\le C$, we have $\|\d\|_2 \le C$ by elliptic estimates.
This implies
$\|a\|_1 \le C$. Applying this, the second equation of (2.5) and 
elliptic estimates, we
deduce $\|\p\|_1\le C$. Higher regularity and estimates 
follow from 
elliptic estimates and imply  the desired compactness.  
\qed
\enddemo

Henceforth we drop the 
subscript 
$1$ in $\SS \W_{1}$, $\R_{1} $ and $\R^0_{1}$.  
It is clear that $\R = \R^0 \slash S^1$, 
where $S^1$ is the group of constant 
gauges.  We deal with $\R^*$ 
 and $\R^0$ separately.

\bigskip
\noindent {\bf The moduli space $\R^*$}
\smallskip
For a given $(a,\p) \in \A_1(Y)\times\Ga_1(Y)$,    let $G_Y=G_{Y,(a, \p)}
: \O_1^0(Y) \to
\O^1_0(Y) \oplus \Ga_0(S)$ be the infinitesimal gauge action operator
at $(a, \p)$, i.e.
$G_Y (f)  = (df,-f \p)$. 
Let $G_Y^* = G_{Y, (a, \p)}^*: \O^1_1(Y)\oplus\Ga_1(S)\to\O_0^0(Y)$ be
the formal adjoint operator
 of $G_Y$ w.r.t. the inner product (2.7). 
We have
$$G_Y^*(b,\psi)=d^*b+ {Im}\langle \p,\psi\rangle.$$

There is a decomposition $\O^1_1(Y) \oplus \Ga_1(S) = 
{ker~ }G_Y^* \oplus {im~ }G_Y$. To 
be more precise, we write 
$\O^1_1(Y) \oplus \Ga_1(S) =
{ker}_1 ~G_Y^* \oplus {im}_1 ~G_Y$.
It follows that the tangent space 
$T_{[a, \p]}\B^*_1(Y)$ of $\B^*_1(Y)$ 
at $[a, \p]$ is represented by 
${ker}_1 ~G_{Y,(a,\p)}^*$ (for any representative 
$(a, \p)$ in $[a, \p]$).
Indeed, the latter gives rise to a vector  
bundle ${Ker}_1 ~G_Y^*  \to \A_1 \times   (\Ga_1(S) - \{0\})$, 
whose quotient bundle ${\Cal K}er_1 ~G_Y^*$ under the action of $\G_2(Y)$ can 
be 
identified 
with  $T\B^*_1$.  
Next note that by 
the gauge invariance of 
the Chern-Simons functional and 
Lemma 2.6, 
we have 
$$
G_Y^* \hbox{\bf sw}_{\lambda, H} (a, \p) =0
$$
in the weak sense for any $(a, \p) \in \A_1(Y) \times 
\Ga_1(Y)$. Hence  
the operator $\hbox{\bf sw}_{\lambda, H}$ defines a section 
$[\hbox{\bf sw}_{\lambda, H}]$ of the quotient bundle ${ ker}_0 ~G_Y^*$,
whose fiber at $[a, \p]$ is represented 
by the weak kernel ${ker}_0 ~G^*_Y$ of $\G_Y^*$ in $\O_0(Y) \oplus
 \Ga_0(S).$ The moduli space $\R_1^*$ 
is precisely the zero locus of this section. 

In the sequel we omit the subscripts $\lambda, H$ in the 
notation $\hbox{\bf sw}_{\lambda, H}$.
Consider the operator $d\hbox{\bf sw}|_{(a, \p)} : \O^1_1(Y) \oplus 
\Ga_1(S) \to \O^1_0(Y)
\oplus \Ga_0(S)$, where $d \hbox{\bf sw }$  means the derivative, i.e. 
the tangent map of the operator $\hbox{\bf sw}$.) By Lemma 2.6, 
it is formally self-adjoint. Assume $[a, \p]  \in \R_1 $. Then 
the gauge invariance of the equation $\hbox{\bf sw} =0$ implies 
$d \hbox{\bf sw} \circ G_Y =0$. It follows that $ G_Y^* \circ d 
\hbox{\bf sw} =0$.
Hence we obtain an operator 
$d \hbox{\bf sw}|_{{ker}_1 ~G^*_{Y, (a, \p)}}: {ker}_1 ~G^*_{Y,(a, 
\p)} \to {ker}_0 ~G^*_{Y, (a, \p)}$ (for 
any representative $(a,\p)$ in $[a,\p]$). It is easy to see that 
this operator is Fredholm of index zero
and represents the linearization of the section 
$[\hbox{\bf sw}]$.
Let it be denoted by $\DD = \DD_{a,\p}$. 
Lemma 2.6 implies that it coincides with the Hessian operator 
of the Chern-Simons functional with respect to the product (2.7).
Note that it extends straightforwardly to 
reducible Seiberg-Witten points. 
\proclaim{Definition 3.3} Let  
$(a, \p)$ be a Seiberg-Witten point, i.e.
a solution of the (perturbed) Seiberg-Witten equation (2.5). 
It is called non-degenerate if
$\DD_{a,\p}$ is onto. The classes $[a,\p]$ or $[a, \p]_0$ are 
called nondegenerate if a representative is nondegenerate.
$($This is independent of the choice of the representative.$)$
\endproclaim

\proclaim{Lemma 3.4} If all elements in $\R^*$ are nondegenerate,
then it is a naturally oriented smooth manifold of dimension zero. $($The 
orientation means that every point in $\R^*$ is assigned a sign.$)$
\endproclaim

\demo{Proof}
By the above discussions, we only need to produce the natural orientation. We can 
use either the degree of the operator $\hbox
{\bf sw}$ or the spectral flow of the operator  $Q$ below as in $\Ta$. 
(They give the same orientation.)
\enddemo

To analyse the operator $\DD$, we introduce another closely 
related formally self-adjoint Fredholm operator $Q$. 
(The Fredholm
property of $\DD$ is also a consequence of the Fredholm
property of $Q$.) First notice the following deformation complex 
$$0\longrightarrow\O^0\mathop{\longrightarrow}^{G_{a,\p}}\O^1
\oplus \Ga 
\mathop{\longrightarrow}^{d\hbox{\bf sw}}
\O^1\oplus \Ga \mathop{\longrightarrow}^{G_{a,\p}^*}
\O^0\longrightarrow 0\tag 3.1$$
where the letters $Y$ and $S$ and the Sobolev subscripts are omitted 
in the notations.
We define $Q=Q_{(a,\p)}:(\O^1_1(Y)\oplus\Ga_1(S))\oplus
\O^0_1(Y)
\to(\O^1_0(Y)\oplus\Ga_0(S))\oplus\O^0_0(Y)$  by the following 
formula:
$$Q=\pmatrix d\hbox{\bf sw} & G\cr
             G^* & 0\cr
\endpmatrix
$$

\proclaim{Lemma 3.5}  Let $(a, \p) \in \SS\W$. 
Then we have
$${ker~}Q\cong \cases\hbox{ker~}\DD\oplus\RR,&
\hbox{if } (a, \p) \hbox{ is  reducible};\cr
{ker~}\DD, &\hbox{if } (a, \p)\hbox{ is irreducible}.\cr\endcases
$$
and
$${coker~}Q\cong \cases{coker~}\DD\oplus\RR,&
\hbox{if } (\p,a) \hbox{ is reducible};\cr
{coker~}\DD, &\hbox{if } (\p,a)\hbox{ is irreducible}.\cr\endcases
$$
\endproclaim

We omit the simple proof. 

\bigskip
\noindent {\bf The moduli space  $\R^0$.}
\smallskip
 The above treatment does not apply to the reducible elements of $\R$, 
 because the tangent bundle
of $\B^*_1(Y)$ does not extend smoothly across the reducibles.
 To analyse the structure of $\R^0$ around reducibles, one can use 
a quotient bundle formulation on the level of the based gauge quotient. 
But we choose a different approach which gives somewhat 
stronger results. Henceforth we make

\proclaim{Assumption 3.6} $Y$ is a rational homology sphere, i.e. its 
first Betti number is zero. \endproclaim

\proclaim{Lemma 3.7} There is a canonical diffeomorphism from 
$\Sigma \equiv (a_0 + {ker}_1~d^*a) \times \Ga_1(S)$ 
to $\B^0_1(Y)$.  In other words, the
former space is a global slice of the action 
of the group $\G_2^0(Y)$ on the space $\A_1(Y) \times \Ga_1(S)$.
\endproclaim 
\demo{Proof} To show that the
natural map from the former space to the latter is one to one, 
consider $b_1, b_2 \in {ker}_1 ~d^*$ and $\p_1, \p_2 \in
 \Ga_1(S)$ such that $(a_0 + b_2,
\p_2) = g^*(a_0 + b_1, \p_1)$ for some gauge $g \in \G^0_2$. 
 Then $d^*(g^{-1}dg)=0$ and $g(y_0)=1$.  Since $Y$ is a rational 
 homology sphere, we have $dg
\equiv 0$, and hence $g \equiv 1$.  The remaining part of
 the proof is obvious.  
 \enddemo
 
  This lemma enables us to reduce the Seiberg-Witten operator to the
said global slice.  But the operator $Q$ is no longer suitable for
 analysing the linearization of the Seiberg-Witten operator.  Instead, 
 we consider the following augmented Seiberg-Witten equation

 $$ \eqalign{ *F_a+df +\langle e_i\cdot\p,\p\rangle
  e^i& = \nabla H(a),\cr \np_a\p + \lambda \p +f \p & =0,\cr}
\tag 3.2 $$
 where $a \in \A_1(Y), \p \in \Ga_1(S)$ and $ f \in \O^0_1(Y)$.

\proclaim{Lemma 3.8} Let  $(a, \p, f)$ be a solution of (3.2 ).
Then $(a, \p)$ satisfies the Seiberg-Witten equation (2.5) and 
$f$ is a constant. Moreover, if $(a, \p)$ is irreducible, then 
$f$ must be  zero. 
\endproclaim
\demo{Proof} Applying $d^*$ to the first equation of (3.2), we 
deduce
$$
d^*df+ f |\p|^2 =0.
$$
The desired conclusion follows.
\qed
\enddemo

We denote the left hand side of (3.2) by
$\hbox{\bf swa}((a, \p, f))$. The linearization
of the restriction of {\bf swa} to $\Sigma$
will be denoted by $\DD_1$. One readily checks that it is a Fredholm
operator of index 1.

\definition{Definition 3.9} Let $(a, \p)$ be a Seiberg-Witten point. 
It is called based-$nondege-$ $nerate$, if $\DD_1$ is onto at 
$(a, \p, f)$, where $f $ is an arbitrary constant 
if $\p =0$ and zero if $\p \not \equiv 0$.
It is easy to see that the based-nondegenerate property is invariant 
under gauge transformations. In particular, this 
definition makes sense for based gauge classes $[a, \p]_0$.
\enddefinition

Let the gauges act on $f$ trivially. The moduli 
space of based gauge classes of solutions of the equation (3.2) 
will be denoted by $\R^0_a$. An element in it is called
based-nondegenerate, if its corresponding element in $\R^0$
is so. As an immediate consequence of the above discussions we obtain 
\proclaim{Lemma 3.10} If all elements of $\R^0_a$ are 
based-nondegenerate, then it is a smooth oriented manifold 
of dimension one. If moreover $\R^0$ is compact, 
then the irreducible part of $\R^0_a$ consists of 
finitely many disjoint circles and its 
reducible part consists of finitely many disjoint
lines. Consequently, the irreducible part of $\R^0$
consists of finitely many disjoint circles
and its reducible part consists of finitely many points.

\endproclaim

Next we give the definition of holonomy perturbations. 
We follow $\BD$  and $\Fl$. Let $D$ denote 
the unit disk in $\Bbb C$. Consider a triple $(y_0, v_0, I)$,
where $y_0\in Y$, $v_0\in T_{y_0}Y$
and $I:D^2\to Y$ is a smooth embedding such that $I(0)=y_o$ and
$dI(T_0D)$ is transversal to $v_0$. 
Fix a point $s_0 \in S^1$. Let $P(y_0,v_0,I)$
be the set of all smooth embeddings  $\ga:S^1\times D\to Y$
such that $\ga(s_0,\t)=I(\t)$ for $\t\in D$ and $\frac{\pa \ga}{\pa s}
(s_0,0)=v_0$. Here $s_0$ is a fixed point in $S^1$. We set
$$P^{(m)}=\cup_{(y_o,v_0,I)}(P(y_0,v_o,I))^m,$$
for $m\in\NN$.  Now we define a map $\ga^h:\A(Y)\times
P^{(m)}\to C^{\infty}(D^2,U(1)^m)$ by
$$\ga^h(a,(\ga^1, \ga^2,\cdots,\ga^m)) (\t)=
(\ga^1_{\t}(a),\ga^2_{\t}(a),\cdot,\ga_{\t}^m(a)),$$
where $\ga^i_{\t}:\A(Y)\to U(1)$ denotes the holonomy map along 
the loop  
$\ga^i(\cdot,\t)$ (at the base point $y_0$).  It is easy to 
see that $\ga^h$ is
gauge invariant. Next  we choose a sequence $\{\e_i\}$ of 
positive numbers
as in $\Fl$  such that
$$C^{\e}(U(1)^m,\RR)=\{u\in  C^{\infty}(U(1)^m,\RR)
: \| v\|_{\e}<\infty\}$$
is complete. Here 
$$\|u\|_{\e}=\sum_{i=0}^{\infty}\e_i\max_{ U(1)^m}|\n^i u|.$$

Now we set
$$\Pi=\cup_{m\in\NN}(P_m\times C^{\e}(U(1)^m,\RR)).$$
This is the parameter space of 
holonomy perturbations. Choose a 
smooth function $\xi$ with support in the interior 
of $D$. For each $\pi=(\ga,u)\in \Pi$, we define
the holonomy perturbation $H_{\pi}:\A(Y)\to\RR$ by
$$H_{\pi}(a)=\int_{D^2}u(\ga^h(a))\xi(\t)d^2\t.$$
It is clear that $H_{\pi}$ extends to $\A_1(Y)$.

\proclaim{Lemma 3.11} For any $\pi=(\ga,u)\in\Pi$, $H_{\pi}$ is a smooth
$\G_{2}(Y)$-invariant function. Moreover, 
the $L^2$-gradient $\n H_{\pi}$ satisfies 
$$\|\n H_{\pi}(a)\|_{L^{\infty}}\le C,$$ 
with $C>0$ independent of $a\in\A$. Similar bounds hold 
for the higher derivatives of $H_{\pi}$. The bounds can 
be made arbitrarily small by choosing $u$ small.
\endproclaim
\demo{Proof} 
For simplicity, we only consider
$m=1$. Set $H=H_{\pi}$. We can write $H(a)=\int_{D^2}u(\ga_{\t}(a))\xi 
d^2\t$. It follows that 
$$dH(a)(b)=\int_{D^2}du_{|_{\ga_{\t}(a)}}(\ga'_{\t}(a)b)\xi d^2\t.$$
Elementary computations lead to
$$\ga'_{\t}(a)b=-\ga_{\t}(a)\int_{\ga_\t}b.$$
We deduce 
$$\eqalign{dH(a)(b)&
=-\int_{D^2}\xi\langle\n u,\ga_{\t}(a)\rangle d^2\t\int_
{\ga_{\t}}b\cr
&=\int_{\ga(S^1\times D^2)}(\xi\circ\ga^{-1})\langle b,\overline{
\langle \n u,\ga_{\t}(a)\rangle}(\ga^{-1})^*(dt)\rangle |\frac{\partial\ga}
{\partial t}|^{-2} (\ga^{-1})^*(dtd^2\t)
\cr
&=\int_Y f (\xi\circ\ga^{-1})\langle b,\overline{
\langle \n u,\ga_{\t}(a)\rangle}(\ga^{-1})^*(dt)\rangle 
\cr
},$$
where
$$ f = |\frac{\partial \ga}{\partial t}|^{-2} |
(\ga^{-1})^*(dt d^2 \t) \slash dvol|.
$$ 
Consequently, $\n H(a)= f (\xi\circ\ga^{-1})\overline{
\langle \n u,\ga_{\t}(a)\rangle}(\ga^{-1})^*(dt)$. 
The desired estimate for $\n H$ follows. The higher order 
derivatives can easily be computed by using the above formula.
\qed
\enddemo

We make 
\proclaim{Assumption 3.12} 
Henceforth we choose 
$H$ in (2.5) and (3.2) to be $H_{\pi}$.
\endproclaim 

We remark in passing that for the purpose  of 
achieving transversality for the moduli 
spaces $\R^*$ and $\R^0$  it is not necessary 
to introduce the holonomy perturbations.
However, they are important for achieving transversality for 
Seiberg-Witten trajectories as will be seen in the next section.

\proclaim{Lemma 3.13} For perturbation
$\pi\in\Pi$ such that $\n^2H_{\pi}$ is small 
enough in $L^{\infty}$-norm (the set of such 
$\pi$ is a nonempty open set), 
there exists a unique reducible element $[a,0] \in\R$.
Equivalently, there is a  
unique 
$a\in\A_1(Y)$ such that
$$\eqalign{ *F_a-\n H_{\pi}(a) & =0,\cr
d^*a & =0.\cr}
\tag 3.3$$
\endproclaim
\demo{Proof} Since $Y$ is a rational homology 
sphere, the 
operator $*d: {ker}~d^* \to {ker}~d^*$
is a bounded isomorphism. Hence the existence 
follows from the implicit function theorem. 
To prove the uniqueness.
consider connections $a$ and $a_1$ satisfying  (3.3).
We set 
$b=a-a_1$ and deduce  
$$*db=\n H_{\pi}(a_1)-\n H_{\pi}(a) \hbox{ and } d^*b=0.$$
By the implicit function theorem, for $\pi$ with 
the property stated in the lemma, $b=0$.
\qed
\enddemo
The unique solution of (3.3) will be denoted by $a_{(h,\pi)}$.

\proclaim{Lemma 3.14} For $\pi\in\Pi$ satisfying the 
condition of Lemma 3.13, let
$\si(\np_{a_{(h,\pi)}})$ be the  set of eigenvalues of 
$\np_{a_{(h,\pi) }}$. 
Assume that $\n^3H_{\pi}$ is small enough in $L^{\infty}$-norm
(the set of such $\pi$ is a nonempty open set). Then for 
$\lmd\in(\RR-\si(\np_{a_{(h,\pi)}})) $, all 
elements of $\R(Y)$ 
are nondegenerate and all elements in $\R^0(Y)$ 
are based-nondegenerate.   
\endproclaim
\demo{Proof} 
We only present the proof for the statement 
concerning the non-degeneracy. The 
based-nondegeneracy can be treated in a similar way.
Consider $(a, \p)\in [a, \p] \in \R(Y)$, we are going to show that $\DD$ 
at 
$(a, \p)$ is onto. By gauge equivariance, we can choose $a = a_{(h, \pi)}$
for the reducible element.
By Lemma 3.5, it suffices to analyse 
the operator $Q$. 
We have 

$$Q(b, \psi, f) =\pmatrix
*db+2\langle e_i\cdot\p,\psi\rangle e^i+ df-\n^2 H(a)b\cr
\np_a\psi+\lmd\psi+f \p+ b\p\cr
d^*b+{Im~}\langle\p,\psi\rangle\cr
\endpmatrix.
$$
Consider an element $(b_1, \psi_1, f_1) \in 
\O^1_0(Y) \oplus \Ga_0(S)\oplus \O_0^0(Y)$ 
satisfying
$$\langle Q(b, \psi, f), (b_1, \psi_1, f_1)\rangle_{L^2}=0\tag3.4$$
for all $(b,\psi, f)\in \O^1_1(Y)\oplus\Ga_1(S) \oplus \O^0_1(Y)$. 
We first derive that $(b_1,  \psi_1, f_1)$ satisfies the adjoint 
equation $Q^* =0$ (hence it satisfies $Q=0$ because $Q^*=Q$)
and is smooth.  
\smallskip
{\bf Case 1} $\p =0$ and $a= a_{(h, \pi)}$.
\smallskip
We have $\np_a \psi_1 + \lambda \psi_1 =0$. By the choice of $\lambda$, we 
conclude that $\psi_1 \equiv 0$.
Now $(b_1, f_1)$ satisfies the following equation
$$\eqalign{*db_1 + df_1 - \n^2H(a)b_1 &=0, \cr
d^*b_1 &= 0. \cr}
\tag3.5
$$
Since $Y$ is a rational homology sphere, the operator $(b_1, f_1) 
\to (*db_1 +df_1, d^*b_1)$ is an isomorphism from $\O^1_1(Y) \oplus   
(\O_1^0(Y))^0$ onto $\O^1_0(Y) \oplus (\O_0^0(Y))^0$, where 
the superscript $0$ means the condition that the average be zero.
As in the proof of Lemma 3.9, we deduce that if $\pi$ has been 
chosen small enough,
$f_1$ must be a constant and $b_1 = 0$.  We conclude that 
${coker~ }Q \cong \RR$. By Lemma 3.5, this implies that $\DD$
is onto.
\smallskip
{\bf Case 2} $\p \not \equiv 0$.
\smallskip
By the unique continuation, the set $U = \{\p \not = 0\}$ 
is an open dense set.  
For $y \in U$, $e_1\cdot\p(y), e_2\cdot\p(y),e_3\cdot\p(y)$ 
and $\p(y)$ span 
$S_y$, 
where
$S_y$ denotes  the fiber of $S$ at $y\in Y$. We deduce that
$\psi_1(y)=0$ for $y\in U$, whence $ \psi_1 \equiv 0$.

Now we easily see that $(b_1, f_1)$ 
satisfies the equation (3.5). Hence $b_1 \equiv 0$ and 
the equation $Q(b_1, \psi_1, f_1) =0$ reduces to 
$f_1 \p =0$.  It follows that $f_1 = 0$ in $U$ and 
consequently $f_1 \equiv 0$. We conclude that 
$Q$ is onto. By Lemma 3.5, $\DD$ is onto.
\qed
\enddemo  

As a consequence of the previous lemmas, we deduce
\proclaim{Proposition 3.15}
Let $\pi$ and $\lambda $ satisfy the same conditions as in 
Lemma 3.13 and Lemma 3.14. Then $\R$ consists of finitely many 
signed points,  the irreducible part of $\R^0$ consists of finitely 
many disjoint oriented circles, and its reducible part is a signed  point. 
\endproclaim

\definition{Definition 3.16} We call 
$\pi$ and $\lambda$ ``good", 
provided that they satisfy 
the conditions of Lemma 3.13 and Lemma 3.14.
The set of good parameters $(\pi, \lambda)$ 
is a nonempty open set.
\enddefinition
  
\head 4. Seiberg-Witten Trajectories: Transversality
\endhead

By Seiberg-Witten trajectories we 
mean solutions of the equation (2.9).
(We shall choose $H = H_{\pi}$.) 
As stated in the introduction, 
our goal is to establish a Morse-Floer 
theory 
for the Chern-Simons functional on the quotient 
space $\B^0_1(Y) = (\A_1(Y) \times \Ga_1(Y)) 
\slash  \G_2^0(Y)$.  
 The union of the critical 
submanifolds of the Chern-Simons functional is 
precisely the moduli space $\R^0$. 
The negative gradient flow
lines of the Chern-Simons functional are given by 
the temporal form of  the Seiberg-Witten trajectories, cf. Section 2.
Setting $A=a(t)+f(\cdot, t)dt$ 
and $\p(t)=\P(\cdot, t)$ in (2.9) we can rewrite it as follows
$$
\eqalign{
\frac{\partial a}{\partial t}-*F_a-d_Yf-\langle e_i\cdot\p,\p\rangle e^i
&=\n H(a), \cr
\frac{\partial \p}{\partial t}+ \np_a\p+\lmd\p+f\p&=0.\cr
}\tag 4.1
$$
(We omit the subscript $\pi$ in $H_{\pi}$.)

We shall use various spaces of local $(l, 4)$ Sobolev class
($L^{l, 4}_{loc}$ class), 
e.g.  $\A_{l, loc} = \A_{l, 4,loc}(X)$, $\Ga^{\pm}_{l,4, loc}$
$ =
\Ga_{l, 4,loc}(W^{\pm}) $, $\O^k_{l, loc} = \O^k_{l,4, loc}(X)$
and $\G_{l, loc}(X) = L^{l, 4}_{loc}(X, S^1)$.

\proclaim{Definition 4.1} We have the following spaces of 
Seiberg-Witten trajectories:

$\N=\{(A, \P) \in \A_{1, loc} \times \Ga_{1, loc}^+: 
(A, \P) \hbox{ solves } (4.1) $
and has finite $($perturbed$)$ 
Seiberg-Witten energy $\}$. 

The corresponding moduli spaces are:
 
$\M=\N
\slash \G_{2, loc}(X)$ and  
$\M^0 = \N \slash \G_{2, loc}^0(X),$  

\noindent
where $\G_{2, loc}^0(X)= \{g \in \G_{2, loc}(X): g((y_0, 0))=1\}.$
$($Recall that $y_0$ is a fixed reference point in $Y$.$)$
\endproclaim

\proclaim{Proposition 4.2} Assume 
that $\lambda$ and $\pi$
are good. Then there are 
positive constants $C$ and $\varepsilon_0$ depending only
on $h$, $\lambda$ and $\pi$ 
with the following properties. For any temporal Seiberg-Witten 
trajectory $u=(A,\P)=(\p,a)$ of 
local $(1,4)$-Sobolev class  and finite energy,
there exist a gauge $g\in \G_2(Y)$ and two smooth solutions
$u_-,u_+$ of $($2.5$)$ such that $g^*u$ is smooth
and the following holds. For all $l$,
$\|g^*u(\cdot,t)-u_{\pm}\|_l<C(l)e^{-C|t|}$ for $|t| \geq T$, where 
$C(l)$ depends only on $h$, $\lambda$, $\pi$, 
$l$ and an upper bound of the energy of $u$, and 
$T >0$ satisfies $E(u, \{|t| >T\}) \leq \varepsilon_0$. 
$(E(u,  \Omega) $ means the energy of $u$ on the domain $\Omega.)$
Moreover, if $u_-$ and $u_+$ are not gauge equivalent, then
$$
E(u) \geq \varepsilon_0.
$$
\endproclaim

The proof will be given in Part II.

\proclaim{Corollary 4.3} We have 
$$
\M= \cup_{\a, \b \in \R} \M(\a, \b), \M^0 =\cup_{\a, \b \in \R}
\M^0(\a, \b),
$$
where $\M(\a, \b) = \N(\a, \b) \slash \G_{2, loc}(X)\hbox{ and }
\M^0(\a, \b)=\N(\a, \b)
 \slash 
\G^0_{2, loc}(X)$ 
and for sets $B_1$ and $ B_2$ of Seiberg-Witten poins, 
$$\eqalign{&\N(B_1, B_2) = \{ u \in \N: \hbox{there is a } 
g \in \G_{2, loc}(X) \hbox{ such that} \cr
&g^*u(\cdot, t) 
\hbox{ converges exponentially to some }  p \in B_1
\hbox{ as } t \to -\infty\cr
&\hbox{and to some } q \in B_2 \hbox{ as } t \to +\infty \}. \cr}\tag 4.2 $$  
\endproclaim

We shall use the moduli spaces $\M^0(\a, \b)$ 
(or the corresponding moduli spaces of trajectories 
in temporal form) to 
construct the boundary operator in our Bott type chain (cochain) complex.

Proposition 4.2 and its corollary suggest that we can work in the set-up 
of exponentially converging trajectories. We introduce various relevant 
spaces. For a positive function $\xi$ on $X$ we consider 
the following $\xi$-weighted $(l, 4)$-Sobolev norms
$$\|u\|_{l,\xi}=(\sum_{k\le l}\int_{X} \xi 
|\n^ku|^{4} dydt)^{\frac{1}{4}}. \tag 4.3
$$
For each  pair of nonnegative numbers
$\d = (\d_-, \d_+)$ we choose a positive smooth function
$\d_{F}$ on $\RR$
such that ${\d}_F(t) =\d_{\pm}|t|$ near ${\pm} \infty$. 
We have the following $\d$-weighted $(l, 4)$-Sobolev
norms ($L^{l, 4}_{\d}$-norms)
$$
\|u\|_{l, \d}=\|u\|_{l, e^{{\d}_F}}.
$$
Let $\O_{l,\d}^k= \O_{l,\d}^k(X)$, 
 $\A_{l, \d}=\A_{l, \d}(X)$ and $\Ga_{l, \d}^{\pm} =\Ga_{l, \d}(W^{\pm})$ 
denote  the completion of the obvious spaces w.r.t. the $L^{l,4}_{\d}$-
norm.

For $u \in \A_1(Y) \times \Ga_1(S)$, let $G_{u}$ denote its isotropy 
 group of gauge actions. (It is trivial if 
$u$ is irreducible and $S^1$ if $u$ is reducible.) 
The isotropy groups are identical for gauge equivalent elements, 
hence $G_{[u]}$ is well-defined.

In the following definition, the notation $\G_{\d}$ should not 
be confused with e.g. $\G_2$. To avoid inconsistence, we 
require that $\d_+$ and $\d_-$ be smaller than 1. 

\proclaim{Definition 4.4} 
For $p, q \in \A_1(Y) \times \Ga_1(S)$ we introduce 
\smallskip
$
L_{\d}(p, q)= 
\{u \in \A_{1, loc} \times \Ga_{1, loc}^+: $
$u-u_0
\in 
\O^1_{1, \d} \times \Ga^+_{1, \d} \hbox{ for some }$ 
$u_0 \in \A_{1, loc} \times \Ga_{1, loc}^+
\hbox{ which near }  \infty (-\infty) $ is t-dependent
and equals $ p (q) \}$.

\smallskip

$ \G_{ \d}(p, q)= \{g\in \G_{2, loc}(X):$ 
$g-g_0\in L^{2, 4}_{\d}(X, \CC) 
\hbox{ for some } g_0 \in \G_{2,  loc}$ which is t-independent 
and  
belongs to $G_{p} (G_q)$ near  
$\infty  (-\infty) \},$ 

$\G_{\d}^0(p, q) = \{g \in \G_{\d}(p, q):
g((y_0, 0)) =1\},$

$\G_{ \d}^I= \{g\in \G_{2, loc}(X):
g-1\in L^{2,4}_{\d}(X, \CC)\}, $

$\G_{\d}^{I, 0} = \{g \in \G_{\d}^I: 
g(y_0, 0)=1\}.$

Note that  
the second and third groups act freely. The first acts freely on 
the irreducible part of $L_{\d}(p,q)$. 
We have the quotients: $\B_{\d} =
\B_{\d}(p, q) = L_{\d}(p, q) \slash \G_{\d}$, 
$\B_{\d}^0 = L_{\d} \slash \G_{\d}^0, $ 
 
and $\B_{ \d}^I = L_{\d} \slash \G_{\d}^I.$ 
\endproclaim

Because $\G_2(Y)$ 
is connected, we obtain  equivalent (in terms of 
suitable gauges) spaces for 
Seiberg-Witten points $p', q'$
which are gauge equivalent to 
$p, q$ respectively.

\proclaim{Definition 4.5} For $\a, \b 
\in \B_1(Y)$ we introduce 

$L_{\d}(\a, \b) = 
\cup_{p \in\a, q \in \b} L_{\d}(p, q)$, 

$\G_{\d}^{\infty} =\{ g \in \G_{2,loc}(X): 
g-g_0 \in L^{2,4}_{\d}(X, \CC)$ for some 
$g_0 \in L^{2,4}_{loc}$ which is $t$-independent 
near $\pm \infty \}$, 

$\G_{\d}^{\infty,0} = \{ g \in \G_{\d}^{\infty}:
g(y_0, 0)=1 \}$.
\endproclaim

\proclaim{Definition 4.6}
For $p, q \in \SS \W$, we set
$\N_{\d}(p,  q)
= L_{\d}(p, q) \cap \N$.
We have the moduli spaces of trajectories
$
\M_{\d}(p, q) = \N_{\d}(p, q) \slash\G_{\d}(p, q)$,
$\M_{ \d}^0(p, q) = \N_{\d}(p, q)
\slash\G_{\d}^0(p, q)$ and
$\M_{\d}^I(p, q) = \N_{\d}(p, q) \slash
\G_{\d}^I.
$
\endproclaim

\proclaim{Definition 4.7}
For $\a, \b \in \R$, we 
set $\N_{\d}(\a, \b) = L_{\d}(\a, \b)\cap \N$. We 
have the moduli spaces:

$\M_{\d}(\a, \b) = \N_{\d}(\a, \b) \slash
\G^{\infty}_{\d}$,
$\M^0_{\d}(\a, \b)= \N_{\d}(\a, \b) \slash
\G^{\infty,0}_{\d}$.
\endproclaim

The irreducible part of e.g. $\M_{\d}$
will be denoted by $\M_{\d}^*$.
The spaces $\M_{\d}(p, q)$ and $\M_{\d}^0(p, q)$
are canonically isomorphic to $\M_{\d}([p], [q])$
and $\M_{\d}^0([p], [q])$ respectively. Moreover,
we have the following easy lemma.

\proclaim{Lemma 4.8} If $(\pi, \lambda)$ is good,
$p, q$ are smooth representatives of $\a$ and $\b$ 
respectively,
and $\d_-, \d_+$ are  
positive and less than the exponent $c$ in Proposition 4.3, then
there are canonical isomorphisms
$$
\M_{\d}(\a, \b) \cong \M(\a, \b), \M^0_{\d}(\a, \b)
\cong \M^0(\a, \b).
$$
We use these ismorphisms to topologize 
$\M(\a, \b)$ and $\M^0(\a, \b)$.
\endproclaim  
 
Thus, $\M(\a, \b), \M_{\d}(\a, \b)$ and $
\M_{\d}(p, q)$ can be viewed as three different models of the 
same space. The same holds for the 
spaces $\M^0(\a, \b)$ etc..
To analyse the structures of 
these moduli spaces, we focus on 
the set-up $\B_{\d}^I(p, q)$ and $\M_{\d}^I(p, q)$.

For $p, q \in \A_1(Y) \times \Ga_1(S)$,  
consider the infinitesimal gauge action 
operator $G_X = G_{X, (A, \P)}: \O^0_{1, \d}
\to \O^1_{0, \d} \oplus \Ga_{0, \d}^+$ at a given 
$(A, \P) \in L_{\d}(p, q)$, $G_X(f) = (df, -f \P)$. Let 
$G^*_X$ be the formal adjoint operator of $G_X$
{\it w.r.t.} the following inner product
$$
\langle(\P_1,A_1),(\P_2,A_2)\rangle_{\d}=
\int_X(Re
\langle\P_1,\P_2\rangle+\langle A_1,A_2\rangle)
e^{\d_F}dydt.
\tag 4.4$$

We have the following elementary lemma, which is 
analogous to Proposition 2a.1 in $\Fl$.

\proclaim{Lemma 4.9} If $\d_-$ and 
$\d_+$ are positive and small emough, then 
$\O^1_{1, \d} \oplus 
\Ga_{1, \d}^+ = \hbox{im}_1 ~G_X \oplus {ker}_1 ~G_X^*$.
Consequently, the tangent space $T_{[A, \P]}\B_{\d}(p, q)$
is represented by ${ ker }~G_X^*$. 
\endproclaim

Now we assume that $\d$ satisfies the condition of Lemma 4.8.
Let ${\Cal U}_{\d} \to \B_{\d}^I$ denote  
the quotient bundle of the trivial bundle $L_{\d}(p, q) \times 
(\O^+_{1, \d} \oplus \Ga_{1, \d}^-) 
\to L_{\d}(p, q)$ 
($\G_{2,\d}$ acts on $\O^+$ trivially and on $\Ga^-$ by
$g^*\Psi=g^{-1}\Psi$).
The (perturbed) 
Seiberg-Witten operator 
$\hbox{\bf SW } = \hbox{ \bf SW}_{\lambda, H}$ (cf. Section 2) induces 
a section $[\hbox{\bf SW}]$ of this bundle. If 
$p, q \in \SS\W$, then 
its zero locus is precisely the moduli 
space $\M_{\d}(p, q)$. 

The linearization of the section $[\hbox{\bf SW}]$ is 
given by the restriction of the operator $d\hbox{\bf SW}|_{(A, \P)}$ to 
${ker}_1G_X$, which will be denoted by $\DD_X=\DD_{X, (A, \P)}$. 
We introduce another closely related operator
$\F_{p, q}:\O^1_{1, \d}\oplus\Ga^+_{1, \d}\to
 \O_{0, \d}^+\oplus\Ga^-_{0, \d}\oplus\O^0_{0, \d}$ 

$$\F_{p, q}=\pmatrix d \hbox{\bf SW}\cr
G^*_X\cr\endpmatrix.
$$

\proclaim{Lemma 4.10} If $p, q \in \SS\W, (A, \P) \in \N_{\d}(p, q)$,
then ${ ker~}\DD_X={ ker~}\F_{p, q}$ and
 ${ coker~}\DD_X={ coker~}\F_{p,q} .$ Consequently, $\DD_X$
is Fredholm iff $\F_{p,q}$ is Fredholm. If they are 
Fredholm, they have the same index.
\endproclaim
\demo{Proof} The kernel 
equality is clear. By the gauge invariance of the 
Seiberg-Witten equation, we have $d \hbox{\bf SW } \circ 
G_X =0$. Applying this and Lemma 4.7,{BBBB Lemma 4.7 is not right BBBB}
 we derive ${ im~ }\F_{p, q} = { im~ }
\DD_X \oplus G^*_X( { im~}_1~ G_X)$. {BBBB I am not sure, pls check BBBB}
 But the second summand equals $\O^0_{1, \d}$.  
\qed\enddemo

\proclaim{Lemma 4.11} Assume that $\d_- $ and $d_+$
are small enough. If $p$ $(q)$
is reducible, we assume in addition 
that $\d_-$ $(\d_+)$ is positive. Then 
$\F_{p, q}$ is Fredholm for all $p, q \in
\SS \W$.
\endproclaim
\demo{Proof}
Let $\O^{k, Y}$ denote the subspace of $\O^k$ consisting of 
forms which do not contain $dt$. Then $\O^+_{1, \d}$ can be identified 
with $\O^{1, Y}_{\d}$. Using this identification, we have 
for a given $(A, \P) = (a + fdt, \p)$
$$\F_{p,q}(b+\ti f dt, \psi)=\dt \pmatrix b \cr \psi \cr \ti f 
\cr \endpmatrix -\pmatrix
 *d_Yb+d_Y \ti f+2\langle e_i\p,\psi\rangle e^i-\n^2 H(a)\cdot b
\cr  -\np_a\psi-\lmd\p-b\p-f\psi- \ti f\p\cr
d_Y^*b - \d_F' \ti f+{Im}\langle\p,\psi\rangle\cr
\endpmatrix.$$
Hence $\F_{p,q}-(\dt-Q+(0, 0, \d_F'))=(0, f, 0)$, where 
$Q$ was defined in Section 3.  Since $f$ decays exponentially, 
its multiplication is a compact operator.  
Now the limits of the operator $Q-(0, 0, \d_F')$ at $\pm \infty$ 
are formally self-adjoint.  Hence we can follow $\Fl$ or $\Bz$  to
show that $\dt+Q -\d_F'$ is Fredholm. Consequently, $\F_{p,q}$ is Fredholm.
\qed\enddemo

Consider a good pair $(\pi_0, \lambda)$. Choose a neighborhood
$\Pi_0$ of $\pi_0$ such that $\Pi_0 \times \{\lambda\}$ 
consists of good pairs and the smallness conditions in Lemma 4.11 
is uniform for all $\pi\in \Pi_0$ (with $\lambda$ fixed).

\proclaim{Proposition 4.12} Assume that $\d_-$
and $\d_+$ are positive, satisfy the above 
smallness condition for all $\pi \in \Pi_0 $ 
and are less than the constant $c$ in Proposition 4.2. Then for generic  
$\pi \in \Pi_0$ the following holds. For all $p, q \in \SS \W$, 
$[ \hbox{\bf SW}]$ is transversal to the zero section, and  hence
$\M_{\d}^I$ is a smooth manifold of dimension 
$ \hbox{ ind } \F_{p, q}$.  Consequently,
$\M_{\d}$ is a smooth manifold of dimension
$$\hbox{ ind }\F_{p, q} -
  \hbox{dim }G_p- \hbox{dim }G_q$$ 
and $\M_{\d}^0$ is a smooth manifold of dimension 
$$\hbox{ ind }\F_{p, q}  -
 \hbox{dim }G_p- \hbox{dim }G_q + 1.$$ 
\endproclaim
\demo{Proof}  First assume that at least one of $p, q$ is irreducible. 
Then all elements of $L_{\d}(p, q)$ are irreducible.
 We extend the  bundle ${\Cal U}_{\d} \to$ 
$\B^I_{\d} \times \Pi_0$  in 
the  trivial way. Then $[\hbox{\bf SW}]$
gives rise to a section of the extended bundle.
By the Sard-Smale theorem, it suffices 
to show that this section is transversal to the zero section, 
which amounts to the surjectivity of the operator 
$\DD_X \oplus d_{\pi}\hbox{\bf SW } $ 
at all $(A, \P)$ which solve  the Seiberg-Witten
equation with parameter $\pi \in \Pi_0$ (and $\lambda$). 
By Lemma 4.8, the latter is equivalent to 
the surjectivity of  the operator $\F_{p, q}\oplus d_{\pi}\hbox{\bf SW}$, 
which in turn follows from the spinor part of the 
transversality argument in $\KM$  (which is similar to the argument  in 
the proof of  Lemma 3.14) and Floer's transversality
argument in $\Fl$  based on holonomy 
perturbations.  This establishes the statement about the space $\M^I_{\d}$.
The statements about $\M^0_{\d}$ and $\M^*_{\d}$
follow via the involved group actions.

 If both $p$ and $q$ are reducible, then they  represent 
the same (unique) reducible element in $\R$. By 
Lemma 2.7, the energy of every Seiberg-Witten 
trajectory equals zero. By the gauge equivariance 
of the operator $\F_{p, q}$, we can use temporal gauges and assume that 
$p=q= (a_{(h, \pi)}, 0)$ and $(A, \P)  \equiv p$. The asserted 
transversality follows from Lemma C.1 in 
Appendix C. Both $\M^I_{\d}$ and $\M^0_{\d}$ are circles in this case.
\qed
\enddemo

Note that by 
gauge equivariance, the transversality property 
is independent of the choice of the representatives
$p, q$ in their gauge classes.  Since $\R$ consists
of finitely many points, for 
generic $\pi$, the transversality property is shared 
by all $p, q$ with $[p], [q] \in \R$.

\definition{Definition 4.13} We shall say that those $\pi$ 
and the corresponding $\lambda$ as described above are {\it generic}. 
\enddefinition 

 By gauge equivariance, $d\SW, G_X$ and $G_X^*$
are equivariant as can easily be verified.   By this and the homotopy 
invariance of Fredholm index, we also have

\proclaim{Lemma 4.14} For 
given $\a , \b \in \R$, $\hbox{ ind }\F_{p, q}$
is independent of the choice of $p \in \a,q \in \b$ and 
$(A, \P) \in L_{\d}(p, q)$. It is also independent of the choice of 
$(\d_-, \d_+)$ $($satisfying the smallness condition 
in Proposition 4.12$)$.
\endproclaim

\noindent {\bf Remark 4.15} 
We have derived the transversality of the 
space $\M^0_{\d}$ in terms of the transversality
of $\M^I_{\d}$.
We present another, more direct argument, which will 
be useful for analysing the moduli
spaces of parameter-dependent trajectories.  
Let $\B^{I, 0}_{\d}$ denote the quotient of $L_{\d}(p, q)$
under the action of the group 
$\G^{I, 0}_{\d}$,   and $\M^{I, 0}_{\d}$
denote the quotient of $\N_{\d}(p, q)$ under the action of 
the same group.  (The element in $\B^{I, 0}_{\d}$
 determined by $(A, \P)$ will be denoted by $[A, \P]_0$.) 
 We need to derive  the transversality
of this later moduli space.  There is a quotient bundle 
$\UU^0_{\d}$ over $\B^{I, 0}_{\d}$ which is analogous
to the bundle $\UU_{\d}$. As before, the operator $\hbox{
\bf SW }$ induces a section $[\hbox{\bf SW}]_0$ of the bundle 
$ \UU^0_{\d}$, whose zero locus is precisely the moduli space
$\M^{I, 0}_{\d}$.

Choose a smooth imaginary valued function $f_0$ on $X$ with 
$f_0(y_0, 0)=\sqrt{-1}$. If both the limits $p$ and 
$q$ are irreducible, we choose $f_0$ to be compactly supported. 
Otherwise, we choose $f_0$ to take the value $\sqrt{-1}$ near 
the infinity which corresponds to the reducible limit. 
We obtain a decomposition
$\Omega^1_{1, \d} \oplus \Ga^+_{1, \d} = \hbox{ im}_1^0~G_X
\oplus \hbox{ ker}_1~G_X^* \oplus 
\hbox{ span }\{G_Xf_0\}$, where $\hbox{im}^0$ means 
the image of those $f$ with $f(y_0, 0)=0$. 
(The third factor is not orthogonal to the first one.)
It follows that the tangent space $T_{[A, \P]_0}\B^{I,0}_{(\d)}$ 
is represented by ${ker}_1~G_X^* \oplus 
\hbox{ span }\{G_Xf_0\}$. Let $P_{f_0}$
denote the orthogonal projection to the orthogonal completement 
of $G_X^*G_Xf_0$ (with respect to the product (4.2)). Then there holds
${ker}_1~G_X^* \oplus \hbox{ span }\{G_Xf_0\}
= {ker}~P_{f_0}G_X^*$. Now the linearization
of the section $[\hbox{\bf SW}]_0$ is given by the restriction of 
the derivative $d\hbox{\bf SW}$
(at $(A, \P)$) to ${ker}~P_{f_0}G_X^*$,
which we denote by $\DD_X^0$. There is an operator $\F^0_{p, q}$
analogous to  $\F_{p, q}$ which is the combination of $d\hbox{\bf SW}$
with $P_{f_0} G^*_X$.  A statement about the relation 
between $\DD_X^0$ and $\F^0_{p, q}$ similar to Lemma 4.10 holds.
In particular, they have the same index. 
The transversality argument in the proof of Proposition 4.11 also 
applies straightforwardly to $\F^0_{p, q}$. Since 
the moduli space $\M^0_{\d}$ is the quotient of $\M^{I, 0}_{\d}$
under the free action of a compact Lie group, the transversality 
of the former follows.  This 
Lie group is trivial if both $p$ and $q$ are irreducible,
and the circle if  at least one of them is reducible.

We note the following relation:
$$ {ind}~\F^0_{p, q} =  {ind}~\F_{p, q} + 1.
\tag 4.3
$$

Finally, we state an important consequence of  transversality.

\proclaim{Lemma 4.16} 
Let $(\pi, \lambda)$ be generic, $p, q \in
\SS\W$ and $u \in \N_{\d}(p, q)$. Choose 
a reference $u_0 \in  L_{\d}(p, q)$. Then 
$d \SW_u$ has a right inverse $Q_u:\O^+_{0,\d}
\oplus \Ga^-_{0,\d} \to \O^1_{1, \d}\oplus \Ga^+
_{1, \d}$ with
$$ \|Q_u\| \leq C,$$
where $C$ depends only on $\|u-u_0\|_{1, \d}$
and $(\pi, \lambda)$. $Q_u$
is equivariant under gauge actions. In particular,
$\|Q_{g^*u}\|\leq \|g\|_{C^1}\|Q_u\|$ for $ g \in \G^{\infty}_{\d}
\cap C^1(X, S^1)$.
\endproclaim
\demo{Proof}  We proceed in  the context of the 
above remark. Let ${\Cal O}_{u}$ denote
the $L^2_{ \d}$-orthogonal completement of
$ ker ~ d\SW_{u}$ in $
{ker~ } P_{f_0}G^*_X$. 
Then the operator $d\SW_{u}|_{{\Cal O}_{
u}}$ is a bounded isomorphism onto
$\O^+_{0, \d}\oplus \Ga^-_{0, \d} $. We define $Q_{u}$ to be its
inverse. The stated norm estimates and gauge equivariance follow readily.
\qed
\enddemo

\head 5. index and orientation
\endhead

Consider a good pair $(\pi, \lambda)$.
Let $O = O_{h, \pi}$ be the unique reducible element in $\R$. 
For $\a\in\R$ we define 
$$\mu(\a)={ind~}\F_{p,q}-1,$$
where $p \in \a, q \in O$.

Note that $\mu $ can easily be extended to all elements of $\B_1(Y)$. It  
depends on $h, \pi$ and $\lambda$. Elementary computation shows 
$\mu(O) =0$. 

\proclaim{Lemma 5.1} For $p, q, r \in \SS\W $ there holds
$$\hbox{ind }\F_{p,r}=\hbox{ind }\F_{p,q}+\hbox{ind }\F_{q,r}-dim~ G_q.$$
An analogous formula holds for the operator $\F^0_{p, r}$. 
\endproclaim
\demo{Proof} 
This is similar to the corresponding index addition formula in 
Floer's theory $\Fl$. Floer's argument can be applied directly. Another 
argument is as follows. Composing with weight multiplication 
operators, we can transform the operators to Sobolev 
spaces without weight. Then the addition formula 
is the  consequence of a linear version of the gluing argument in Part II. 
The term ${ dim~ }G_q$ arises 
because of  the ``jumping" across the kernel 
of the operator $d^* + { Im}\langle \p, \rangle$
which is caused by the operator $(0, 0, \d_F')$.
\enddemo

\proclaim{Corollary 5.2} There holds 
$$
\hbox{ind }\F_{p, q} = \mu([p]) - \mu([q])+ \hbox{dim }G_q.
$$
Consequently, if $(\pi, \lambda)$ is 
a generic pair, then we have $\hbox{ dim } \M_{\d}(p, q) =
\mu([p]) -\mu([q])  - 
\hbox{dim }G_p,$
$\hbox{ dim } \M^0_{\d}(p, q) = \mu([p]) - \mu([q])+1 
 - \hbox{dim }G_p.$ 
\endproclaim
 
Next we study the orientation of the moduli spaces of 
Seiberg-Witten trajectories. 

\proclaim{Proposition 5.3}Assume that $(\pi, \lambda)$ 
is generic. Then $\M_{\d}^I(p,q), 
\M_{\d}(p, q)$ and $\M_{\d}^0(p, q)$
are orientable.  Indeed, their orientations are canonically determined 
after some choices are made, which will be 
given in the proof below. $($Consequently,
$\M_{\d}([p], [q]),$ $ \M_{\d}^0([p], 
[q]), \M([p], [q])$ and $\M^0([p], [q])$
are orientable.$)$ Moreover, the orientations 
are consistent with the gluing construction 
used in the proof below, namely the orientation 
of $\M^I_{\d}(p, q)$ is the same as the product 
orientation induced from gluing $\M^I_{\d}(p, r)$
to $\M^I_{\d}(r, q)$.
\endproclaim
\demo{Proof} 
Without loss of generality,
assume that $\d_+$ and $\d_-$ are equal and 
positive. The operator $\F_{p, q}$  induces a section of Fredholm 
operators over $\B^I_{\d}$. 
Let $det(p,q)=det~\F_{p,q}$ be its determinant line bundle. Indeed, the 
operator section $\F_{p, q}$ can be deformed through
Fredholm operator sections to the operator section
$(d^+ + d^*_{\d}, D_A)$, whose determinant line bundle 
is trivial. Hence $det(p, q)$ is 
trivial. An orientation of the vector space    $H^0_{\d}
\oplus H^1_{\d} \oplus H^+_{\d}$ (the homology of the complex 
associated with the operator $d^+ +d^*_{\d}$) then determines 
an orientation of $det(p, q)$, cf. $\WI$.

To obtain consistent orientations, we
choose a smooth irreducible pair $p_0=(a_0, \p_0)$  
and consider the space $L_{\d}(p_0, p)$ and its 
quotient $\B^I_{\d}(p_0, p)$ for $p \in \SS\W$.
By the proof of Lemma 4.11, the operator section $\F_{p_0, p}$ 
is a Fredholm operator section. By the above argument,
its determinant line $det(p_0, p)$ is trivial. We fix an orientation 
(a trivialization) for it. For $p, q \in \SS\W$
we construct an embedding by a simple gluing process 
$$ L_{\d}(p_0, p) \times L_{\d}(p, q)
\to L_{\d}(p_0, q). 
$$ 
(Compare $\Fu$.) 
On the other hand, we choose reference elements 
$u_0 \in L_{\d}(p_0, p)$ and  
$u_1 \in L_{\d}(p, q)$. 
Then it is easy to show that $ u_0, u_1   + 
(\hbox{ker } d^*_{\d} \times \Ga_{1, \d}^+)$ 
are global slices in $L_{\d}(p_0, p)$ and 
$L_{\d}(p, q) $ 
for the action of the groups $\G^I$ respectively.
Using them and the above embedding we obtain an embedding $\Theta$:
$$\B_{\d}^I(p_0, p)\times \B_{\d}^I(p, q) 
\to\B_1(p_0, q).$$
We have the projections $\pi_{0}$, 
and $\pi_1$ 
of the above product to its factors.  In addition let
$\pi_p$ be its projection to $p$. 
Now the index addition formula (Lemma 5.1) leads 
to an addition formula for the index 
bundle, which in turn implies a product 
formula for the determinant line bundle. 
We apply the last formula to the present situation to deduce 
$$\pi^*_0det(p_0, p)\otimes \pi_1^* \det(p, q) 
\otimes 
\pi^*_pl_p 
\cong det(p_0, q)|_{\hbox{ im }\Theta},$$
where $l_p$ is the dual of the 
kernel of the operator $d^* + \langle \p, \cdot \rangle$
at $p$. 
We choose an orientation for $l_O$ (note that 
the $l_p$'s are canonically equivalent
to each other for $p \in O$). Then the above isomorphism determines 
an orientation of $det(p, q)$,  which gives rise 
to an orientation of $\M^I_{\d}(p, q)$. The desired consistency 
follows from the construction. 

If both $p$ and $q$ are reducible, then 
$\M^0_{\d}$ is a circle generated by gauge actions,
and hence inherits a canonical orientation from the actions.  
Otherwise, the moduli  space $\M_{\d}$ 
is the quotient of $\M^I_{\d}$ by a free $S^1$ action 
if one of $p, q$ is reducible, and 
equals $\M^I_{\d}$ if neither is reducible. 
Hence the orientation of $\M^I_{\d}$
induces an orientation of $\M_{\d}$.
On the other hand, $\M_{\d}$ is the quotient 
of $\M^0_{\d}$ by another free $S^1$ action, 
hence we arrive at an orientation of $\M^0_{\d}$.
\qed
\enddemo

\proclaim{Remark 5.4} As the transversality of 
$\M^0_{\d}(p, q)$ can be derived directly through
the study of the operator $\F^0_{p, q}$, so can 
its orientability. The orientability part of the above 
arguments can be modified to 
suit the situation of $\M^0{(p, q)}$. But we 
need a different arrangement for consistent 
orientation here, which will be discussed in Section 7. 
\endproclaim

\head 6. The temporal model 
and compactification
\endhead

A temporal Seiberg-Witten trajectory is a solution $(A, \P)$
of (4.1) which is temporal (in temporal form,or ``in
temporal gauge"), i.e. $A = a + fdt$ 
with $f \equiv 0$. For  $ B_1, B_2 \subset \SS\W$ we set 
$$
\N_T(B_1, B_2) = \{u \in \N(B_1, B_2): u \hbox{ is 
temporal} \}, 
$$
$$
\N_{T,{\d}}(B_1, B_2) = \{u \in \N_{\d}(B_1, B_2):
u \hbox{ is temporal}\}.
$$
For $\a  \in \R$ let $S_{\a}$ denote its 
lift to $\R^0$.  For $\a, \b \in \R$ we  set
$$
\M_T( S_{\a}, S_{\b}) = \N_T(S_{\a}, S_{\b})
\slash \G_2^0(Y),
$$
$$
\M_{T,{\d}}(S_{\a}, S_{\b}) =
\N_{T,{\d}}(S_{\a}, S_{\b}) \slash \G_2^0(Y).
$$

By Proposition 4.2, for good parameters and small $\d$,
the first two spaces are identical, and the last two are so, too.
 We shall only consider good parameters and small $\d$. Note that 
we can also assume that the space $\M(p, q)$ is identical to 
the space $\M_{\d}(p, q)$.

\proclaim{Lemma 6.1} 
Assume that $(\pi, \lambda)$ is generic. Then there are canonical 
diffeomorphisms  from $\M^0(p, q)$ and $\M^0(\a, \b)$ to 
$\M_T(S_{\a}, S_{\b})$ for any $p \in \a, q \in \b$, 
where $\a, \b \in \R$. (Indeed, these  diffeomorphisms prove the manifold 
structure of $\M_T(S_{\a}, S_{\b})$.)
\endproclaim
\demo{Proof} 
Let $L_{{\d}, T}(\a, \b)$ denote the temporal part of $L_{\d}(\a,
\b)$, which is a smooth submanifold as can easily be seen. 
We have the following temporal transformation 
$T_G$: for $u =(a + f  dt,\P)$, set
$$g_T(u) = e^{-\int^t_0f}
$$
and $T_G (u) = g_T(u)^*u$. Clearly, $T_G: \M^0(\a, \b)
\to  L_{\d, T}(\a, \b) \slash \G^0_2(Y)$
is a smooth map with image $\M_T(S_{\a}, S_{
\b})$. One readily checks that this map is an embedding. 
\qed
\enddemo

We need to compactify our moduli spaces of trajectories. 
First we introduce some terminology.

\definition {Definition 6.2} 1. Let $p, q \in {\Cal S}{\Cal W}$.
A  {\it  k-trajectory} $u = (u_m)_{1\leq m \leq k} 
$ from $p$ to $q$ with consecutive 
{\it junctures} $p_0=p, ..., p_k=q \in \SS\W$ is an element in 
$\N(p_0, p_1) \times \N(p_{k-1}, p_k)$ with $u_m \in \N(p_{m-1}, p_m)$.
$u^m$ is called the {\it m-th portion} of $u$. 
$p, q$ are called the {\it endpoints} of $u$ at 
$+\infty$ and $-\infty$ respectively.  $u$ is called 
{\it proper}, if its junctures belong to distinct gauge 
classes with respect to the full gauge group. 
\enddefinition

\definition{Definition 6.3} For $u \in \N(p, q)$,
let $p', q'$ be the endpoints of $T_G(u)$ at
$-\infty$ and $+\infty$ respectively. We set
$\pi^-(u) =[p']_0, \pi_+(u) = [q']_0$. $\pi_+$
and $\pi_-$ are called the {\it temporal endpoint maps}.
\enddefinition

\definition{Definition 6.4} 
A {\it k}-trajectory $(u_m)$ is called
{\it consistent}, if $\pi_+(u_m)=\pi_-(u_{m+1})$
for all $1 \leq m \leq k-1$. For distinct
$p_0, ..., p_k \in \SS\W$, let  $\N(p_0, ..., p_k)$
denote the space of consistent (and proper)
{\it k}-trajectories with junctures $p_0, ..., p_k$. 
For $p, q \in \SS\W_0, p \not = q$, let  
$ \N(p, q)^k $ denote the space of proper and consistent
{\it k}-trajectories from $p$ to $q$. We set 
$\hat \N(p, q)= \cup_k \N(p, q)^k$.  Note that 
the temporal endpoints maps $\pi_+$ and $\pi_-$ naturally extend to 
consistent $k$-trajectories. 
\enddefinition

\definition{Definition 6.5} Let e.g. $\hat \N_T(p, q)$ denote 
the subspace of $\hat \N_T(p, q)$ consisting 
of temporal {\it k}-trajectories. For $\a, \b \in \R$ we  
set $\hat \N_T(S_{\a}, S_{\b}) = \cup_{p \in \ga \in S_{\a}, q \in
\ga' \in S_{\b}}$ $ \hat \N_T(p, q)$ and 
$$\hat \M_T(S_{\a}, S_{\b})=
\hat \N_T(S_{\a} , S_{\b}) \slash \G_2^0(Y).$$ (The action occurs on
each portion of $k$-trajectories.)
We also have subspaces $\hat \M_T(S_{\a}, S_{\b})_k$ 
of $\hat \M_T(S_{\a}, S_{\b})$, which  provide it 
with a natural stratification.  Furthermore, 
we have subspaces $\hat \M_T(S_{\a_0},..., S_{\a_k})$
for distinct $\a_0,...,\a_k \in \R$.
\enddefinition

\definition{Definition 6.6} 
For each $\a \in \R$ we choose an element $p_{\a} \in 
\a$. We fix this choice henceforth and denote 
the set of these elements by $\SS \W_0$. 
Let $\N(p, q;\SS\W_0)^k$ denote the subspace of 
$\N(p, q)$ consisting of {\it k}-trajectories 
whose junctures belong to $\SS\W_0$. Similarly, we 
have $\N(p, q;\SS\W_0)$. 

For distinct $p_0, .., p_k \in \SS\W$, we set 
$$ \M^0(p_0,..., p_k)= \N(p_0,..., p_k) \slash 
(\G_{\d}^0(p_{0}, p_1) \times...\times \G_{\d}^0(p_{k-1}, p_k)).$$
For $p, q \in \SS\W_0$, let $ \M^0(p, q)^k$ denote the union of 
all $\M^0(p_0,..., p_k)$ with distinct $p_0=p,..., p_k=q \in \SS\W_0$, 
and let $\hat \M^0(p, q)$ denote the union of $\hat \M^0(p, q)^k$
over all possible $k$. By the definition, $\hat \M^0(p, q)$
has a natural stratification.
\enddefinition

\proclaim{Lemma 6.7} Assume that $(\pi, \lambda)$ 
is generic.  Then $\M^0(p_0,...,p_k)$ and $\M_T(S_{\a_0},
...,$ $S_{\a_k})$ are canonically diffeomorphic smooth manifolds, 
where $\a_0,...,\a_k \in \R$ are distinct and 
$p_i \in \a_i$.
\endproclaim

\proclaim{Proposition 6.8}
Assume that  $\pi, \lambda$ are good. Then there is a $k_0$ 
depending only on $h, \pi $ and $\lambda$ 
such that that the following holds. Let $u^j$ be a sequence of 
temporal Seiberg-Witten trajectories
with uniformly bounded energy. Then 
there is a sequence of gauges 
$g_j \in \G_2^0(Y)$ with the following properties:

1. Set $\tilde u^j =
g_j^*u_j$. Then each $\tilde u^j$ is smooth and 
and converges exponentially at time infinities as described 
in Proposition 4.2.  Let $p_j$ and 
$q_j$ denote its limits at $+\infty$ and $-\infty$ respectively. 

2.  After passing to a subsequence of 
$u^j$, $p_j$ converge to a $p$, and $q_j$ converge 
to a $q$. Moreover, under the assumption that $p$ differs 
from $q$, $\tilde u^j$ converge to a proper 
$k$-trajectory $u=(u_m) \in \hat \N_T(p, q)^k$  with $k
\leq k_0$ in the following sense:

i)$E(\tilde u^j) \to E(u) \equiv\sum_m E(u_m)$,

ii) there are decompositions $\RR = \cup_{1\leq m \leq k} [t_m^j,
t_{m+1}^j]$ with $t_m^j \in [-\infty, +\infty], t_{m+1}^j
>t_m^j$ and $t_{m+1}^j-t_m^j 
\to +\infty$ as $j \to \infty$, such that for each $l \geq 0$, 
$$
sup_{Y \times [t_m^j-T_m^j, t_{m+1}^j-T_m^j]
} e^{\delta_F}|\n^lu^j(\cdot 
+T_m^j, \cdot)-
\n^l u_m| \to 0 \tag 6.1
$$
as $j \to \infty$, where $T_m^j$  are constants with the property 
$t_m^j-T_m^j \to -\infty$ and $t_{m+1}^j-T_m^j
\to +\infty$ as $j \to \infty$, and the $\delta_+$, $\delta_-$ in
$\delta$ are positive and less than the exponent $C$ in Proposition 4.2.

$($Note that 
i) is actually implied by ii).$)$\endproclaim

In short, the above proposition says that 
each $\tilde u^j$ splits into 
$k$ portions to yield $k$ new sequences,
which converge in exponentially weighted norms 
after suitable time translation adjustments.
The fact that the total limit is a $k$-trajectory implies in particular 
that adjacent portions converge to trajectories whose 
endpoints match. This is an important property.

\proclaim{Corollary}
Let $u^j \in \N(p, q)$ with $p, q \in \SS \W_0,  p \not 
= q$. Then there are gauges $g_j \in\G_{2, 
\delta}^0(p, q)$ such that $g_j^* u^j$
are smooth and after passing to a subsequence, they converge 
to a proper and consistent $k$-trajectory $u \in 
\hat \N(p, q; \SS\W_0)$ with $k\leq k_0$ in the sense as 
described in the proposition.
\endproclaim

The proof of this proposition will be given in Part II. 

We extend the concept of convergence of trajectories 
to $k$-trajectories to convergence of 
$k$-trajectories to $k'$-trajectories in the obvious
way. Then Proposition 6.8 and 
its corollary extend straightforwardly to $k$-trajectories. 

The real line $\RR$ acts on trajectories in terms of 
 the time translation.  We define the time translation 
 action on $k$-trajectories to be the separate time translation 
 action on each portion of $k$-trajectories. It 
 gives rise to an action of $\RR^k$. Let the underline denote quotient
under the time translation action in the 
context of temporal trajectories, e.g. ${\underline \M}_T(S_{\a},
S_{\b})= \M_T(S_{\a}, S_{\b}) \slash\RR$.

\proclaim{Proposition 6.9}
Assme the same as in Proposition 6.8. 
Then the spaces ${\hat {\underline  \M}}_T(S_{\a},$  
$ S_{\b})$ are compact, where the topology is given by the 
convergence concept in Proposition 6.4 and its corollary.
Similarly, for $p, q \in \SS \W_0$, the quotient of ${\underline 
{\hat \M^0}}(p, q)$ under the time translation action is compact.  
\endproclaim

Consider $p \in \ga \in S_{\a}, q \in \ga' \in S_{\b}$. The 
spaces $\hat \M^0(p, q)$ (or $\M^0(p,q)$) and ${\hat 
\M}_T(S_{\a}, S_{\b})$ (or $\M_T(S_{\a}, S_{\b})$) are 
isomorphic. But their quotients under the time translation 
action are not isomorphic. For our purpose, the time translation action 
on the temporal model is more suitable. 
For this reason, we consider the action of $\RR$ on $\N(p, q)$ 
induced from the time translation action on $\N_T(S_{\a}, S_{\b})$, 
which we call the ``twisted time translation".

\definition{Definition 6.10}
For $R \in \RR$, let $\tau_R$ denote the time translation 
by $R$, i.e.
$$
\tau_R(u)(y, t) = u(y, t-R).
$$
The twisted time translation $ T_R$ by $R$ 
is defined as follows. Let $u \in \N(p, q)$. Then $ T_R u  
= (g_T(u)^{-1})^*(\tau_R(T_G(u)))$. (See the proof of Lemma 6.1 for 
$g_T$ and $T_G$.) The twisted time translation acts on each portion 
of $k$-trajectories separately, giving rise to a $\RR^k$ action.
We use the underline to denote the quotient under the twisted 
time translation action.  
\enddefinition

\proclaim{Lemma 6.11} Assume the same as in Proposition 6.8. Then 
$\underline {\M^0}(p_0,...,p_k)$ and $\underline {\M}_T(S_{\a_0},
...,S_{\a_k})$ are canonically diffeomorphic 
smooth manifolds, where $p_i, \a_i$ are the same as in Lemma 6.7. 
\endproclaim

The following proposition is a consequence of Proposition
6.4. We assume $p , q \in \SS\W_0, p \not = q$.

\proclaim{Proposition 6.12} Assume the 
same as in Proposition 6.8. Then 
the spaces $\underline {\hat \M^0}(p, q)$ are compact. 
More precisely, consider e.g. a sequence $[u^j]_0 \in 
\M^0(p, q)$.  There is a  sequence $v^j \in \N(p, q)$ such that 
\roster
\item"1" $T_{c_j}v^j \in [u^j]_0$ for 
some $c_j$.
\item"2" after passing to a subsequence, 
$v^j$ converges to a proper and 
consistent $k$-trajectory 
$v=(v^m)  \in \hat \N(p, q; \SS\W_0)$. 
\endroster
\endproclaim

Notice that the compact space ${\underline {\hat \M}}_T(
S_{\a}, S_{\b})$ contains ${\underline \M}_T(S_{\a},
S_{\b})$ as a subspace. The compactification of ${\underline \M}_T(S_{\a},
S_{\b})$ is given by its closure 
${\underline {\bar \M}}_T(S_{\a}, S_{\b})$. 

\proclaim{Proposition 6.13} In a generic situation,
we have ${\underline {\bar \M}}_T(S_{\a},
S_{\b}) = {\underline {\hat \M}}_T(S_{\a}, S_{\b})$.
Moreover, the following hold:

(1) ${\underline {\hat M}}_T(S_{\a},S_{\b})$ has the structure of 
$d$-dimensional smooth orientable  manifolds with corners (i.e. 
modeled on the first quadrant of $\RR^d$), where 
$d = \mu(\a)-\mu(\a)-\hbox{dim }G_p +1$ with $[p] \in \a, [q]\in \b$.

(2) This structure is compatible with the stratification 
${\underline {\hat \M}}_T(S_{\a}, 
S_{\b})\! = \!\cup_k {\underline {\hat \M}}_T(S_{\a},$
$S_{\b})^k$, i.e. the interior of the $k$-dimensional edge of 
${\underline {\hat\M}}_T(S_{\a}, S_{\b})$ 
is exactly ${\underline {\hat \M}}_T(S_{\a},$ $S_{\b})^k$.

(3) The temporal endpoint maps $\pi_-:
{\underline {\hat \M}}_T(S_{\a},
S_{\b}) \to  S_{\a}$ and $\pi_+: {\underline {\hat \M}}_T(S_{\a},
S_{\b})$ $ \to S_{\b}$ are smooth fibrations. $($They 
are naturally induced from the previous 
temporal endpoint maps.$)$

In particular, we have 
$$\eqalign{
\pa \underline {\hat \M}_T(S_{\a}, S_{\b})&=
\cup_{\mu(\a) > \mu(\gamma) > \mu(\b)}
\underline {\hat \M}_T(S_{\a}, S_{\ga},
S_{\b})\cr
&=\cup_{\mu(\a) > \mu(\gamma) > \mu(\b)}
\underline {\hat \M}_T(S_{\a}, S_{\ga}) 
\times_{S_{\ga}} \underline {\hat \M}_T(
S_{\ga}, S_{\b}),\cr} \tag 6.2 
$$
where the fiber product space $\underline {\hat 
\M}_T(S_{\a}, S_{\ga}) \times_{S_{\ga}}
\underline {\hat \M}_T(S_{\ga}, S_{\b})$
is defined to be 
$$\{(u, v)\in \underline {\hat \M}_T(S_{\a}, S_{\ga}) \times
\underline {\hat \M}_T(S_{\ga}, S_{\b}):
\pi_+(u) = \pi_-(v)\}. \tag 6.3
$$
\endproclaim

The following is an equivalent formulation for 
Proposition 6.13. (We only give part of the statements.)

\proclaim{Proposition 6.14} In a generic situation, we have 
${\underline {\bar M^0}}(p, q) = \underline {\hat M^0}(p, q)$.
Moreover, ${\underline {\hat \M^0}}(p,q)$ has the structure of 
$d$-dimensional manifolds with corners which is 
compatible with its natural stratification.
\endproclaim 

The said equivalence means the following:
\proclaim{Lemma 6.15}
In a generic situation, ${\underline {\bar \M^0}}(p, q)$
is canonically diffeomorphic to ${\underline {\hat M}}_T( S_{\a}, S_{\b})$.
\endproclaim

This lemma is an easy consequence of the temporal transformation once 
Proposition 6.14  has been established. Hence Proposition 6.13  
can be derived as a corollary 
of Proposition 6.14. It may be possible to 
prove it without appealing to Proposition 
6.14, but that seems rather cumbersome. Note that Proposition 6.14  alone 
would suffice for our purpose, but Proposition 6.13 is important 
from a conceptual point of view.

Proposition 6.14 is a consequence of the compactness result Proposition 6.7
and a result on structures near infinity which we present now.

We only consider generic situations. Let $u$ be a Seiberg-Witten 
trajectory of finite energy. We define $\rho_+(u)$ and $\rho_-(u)$ 
by the following equations
$$
E(u, \{t \geq \rho_+\})= \varepsilon,
E(u, \{t \leq \rho_-\}) = \varepsilon, 
$$
where $\varepsilon$ is given in Proposition 4.2.
Let $<u>$ denote the set of trajectories which 
are gotten from $u$ by a twisted time translation.
Let $u^*_+(u) $ denote the element in $<u>$ whose $\rho_+$ 
value equals zero, and $u^*_-(u)$ denote that with the  
$\rho_-$ value equaling zero. By a real analyticity argument, 
one readily shows that $u^*_+$ and $u^*_-$ give rise to transversal 
slices for the twisted time translation 
action. We set $\N^*(p, q) = \{ u \in \N(p, q): \rho_-(u) = 0\}$.
Similarly, we have $\N^*(p, q)^k$ etc..

For each $u$, there are unique numbers $R_+(u)$ and $R_-(u)$ such that 
$u = T_{R_+(u)}(u_+^*$ $(u))$ and $u = T_{R_-(u)}(u_-^*(u))$. 
In particular, $\N(p, q)^k$ is isomorphic to $\N^*(p, q)_k \times  \RR^k$.

\proclaim{Proposition 6.16}
There are neighborhoods $U$, $\hat U$ of $\underline {\M^0}(p, q)^k$
in $\underline {\M^0}(p, q)^k\times [0, \infty)^{k-1}$ and 
$\underline {\hat \M}(p, q)$ respectively, and a
homeomorphism $\hbox{\bf F}:U \to \ti U$ such that
the restriction of $\hbox{ \bf F }$ to $U_0=U \cap (\underline
{\M^0}(p, q)^k \times (0, \infty)^{k-1})$
is a diffeomorphism. Moreover, the following hold:

1)
For each compact set $K$ in $\underline {\M^0}(p, q)^k$, there is
a positive number $r_0$ such that $K \times [0, r_0]^{k-1} \subset U$.

2)
For $1 \leq j \leq k-1$,
the restriction of $\hbox{\bf F}$ to the $j$-th boundary
stratum of $U$ is a diffeomorphism onto the $j$-th boundary 
stratum $\ti U \cap\underline {\M^0}(p, q)^{j+1}$
of $U  \cap \underline
{\M^0}(p, q)$. Here e.g. the first
boundary stratum of $\underline
{\M^0}(p, q)^k \times
[0, \infty)^{k-1}$ is
$\underline {\M^0}(p, q)^k  \times (\{0\} \times  (0, \infty)^{k-2}
\cup (0, \infty) \times \{0\} \times (0, \infty)^{k-2}
...\cup (0, \infty)^{k-2} \times \{0\}).$

$\hbox{\bf F}$ defines the structure of
smooth manifolds with corners for $\underline {\hat \M^0}(p, q)$
stated in Proposition 6.13.
\endproclaim

The proof of this proposition will be given in Part II.

\head 7. Bott-type and stable equivariant homology
\endhead

We first introduce a few orientation conventions. We follow those used in 
$\FU$.  For an oriented smooth manifold with 
corners $\X$, its boundary is oriented in 
such a way that $span\{n_{\pa \X}\} \oplus
T\pa \X = TX|_{\pa \X}$ as oriented vector bundles (away from the corners 
of $\pa \X$), where $n_{\pa \X}$ is an inward normal field of the 
boundary. Given transversal smooth maps 
$F_1: \X_1 \to S$ and $F_2: \X_2 \to S $ from 
two oriented smooth manifolds with corners into 
an oriented smooth manifold $S$, the  fiber 
product $\X_1 \times_S \X_2 = (F_1 \times F_2)^{-1}
(Diag (S \times S))$ ($Diag$ means {\it diagonal})
has a canonical orientation such that $T(\X_1 \times_S 
\X_2) \oplus N = (-1)^{dim S \cdot dim \X_2} TX_1 \oplus TX_2$
as oriented bundles, where $N$ denotes the 
oriented bundle $(d (F_1 \times F_2))^{-1}((TS \oplus \{0\})|_{Diag(
S \times S)})$. The following lemma is due to Fukaya and easy to verify. 

\proclaim{Lemma 7.1} 
There hold
$$
\pa(\X_1 \times_S \X_2)= \pa \X_1
\times_S \X_2 + (-1)^{ dim  \X_1
+  dim  S} \X_1 \times_S \pa \X_2, \tag 7.1
$$
and
$$ (X_1 \times_{S} \X_2) \times_{S'} \X_3 =
X_1 \times_S (\X_2 \times_{S'} \X_3). \tag 7.2
$$
\endproclaim 

We shall use the natural orientations of $S_{\a}$ 
($\a \in \R$)  provided by Lemma 3.10. 
(We can also use any other orientations.)
By Proposition 6.13, the boundary of 
$\underline {\hat \M}_T(S_{\a}, S_{\b})$ is 
a union of fiber products. We need to arrange the orientation 
of these spaces so that a suitable consistency holds in regard 
of the natural orientation of fiber products as defined 
above and boundary orientations. Indeed, we have 

\proclaim{Lemma 7.2} We can choose the orientation of 
$\underline {\hat \M}_T(S_{\a}, S_{\b})$ such that
$$
\pa \underline {\hat \M}_T(S_{\a}, S_{\b}) = 
(-1)^{\mu(\a) + dim S_{\a}} \bigcup_{\mu(\a)  > \mu(\ga) 
> \mu(\b)} \underline {\hat \M}_T(S_{\a}, 
S_{\ga}) \times_{S_{\ga}} \underline {\hat \M}_T(S_{\ga},
S_{\b})
\tag 7.3
$$
as oriented manifolds.
\endproclaim  

This lemma is analogous to Sublemma 1.20 in $\FU$ and can be  proven 
by the same arguments as there.

\bigskip

\noindent {\bf Bott-type Seiberg-Witten Floer homology }

\bigskip

For a topological space $\X$, let $ C_j(\X)$ denote
the free abelian group of singular
$j$-chains in $\X$ with coefficient $\ZZ$. Let $A_j(\X)$ be its subgroup 
generated by elements of the form $(\Delta_j, f) - (\Delta_j,
f')$, where $\Delta_j$ denotes the standard Euclidean 
$j$-simplex and $f$, $f'$ are homotopic continuous 
maps from $\Delta_j$ to $\X$. Let $Br_j(\X)$ denote its subgroup generated 
by elements of the form $(\Delta_j, f)-\sigma$, where 
$\sigma$ is obtained from $(\Delta_j, f)$ by a baricentric subdivision.
We set
$$
\ti C_j(\X) =  C_j(\X) \slash (A_j(\X)+
Br_j(\X))
$$
and define $ \ti C^j(\X)$ to be the dual of $\ti C_j(\X)$, i.e. the free 
abelian group of homomorphisms from $\ti C_j(\X)$ to $\ZZ$. 
For $\sigma \in C_j(\X)$, its equivalence class in $\ti C_j(\X)$ 
will be denoted by $<\sigma>$. Let $\partial_{O}$ denote the ordinary 
boundary operator on $C_j(\X)$. 

\definition{Remark 7.3} For $\sigma= (\Delta_j, f) \in C_j(\X)$,
the map $( \pa \Delta_j, f|_{\pa \Delta_j})$
from the oriented boundary induces a chain in 
a natural way. By our convention for boundary 
orientation, this chain equals $- \pa \sigma$.
\enddefinition

Next we set $S_i=\cup \{S_{\a}: \a \in\R, \mu(\a)=i\}$ 
and $$\ti C_k=\oplus_{i+j=k}\ti C_j(S_i),
\ti C^k=\oplus_{i+j=k}\ti C^j(S_i),$$

$$\ti C_* = \oplus_k \ti C_k, \ti C^* = \oplus_k \ti 
C^k.$$

We proceed to define a boundary operator
$\ti \pa: \ti C_k\to \ti C_{k-1}$ 
with the dual (coboundary operator) 
$\ti \pa^*: \ti C^{k-1} \to \ti C^k$ for each $k$. 
First, for each pair $\a, \b \in 
\R$ with $\mu(\a) > \mu(\b)$ we define a boundary operator 
$\pa_{\a, \b}: C_k \to C_{k-1}$, where
$$C_k = \oplus_{i+j=k}C_j(S_i).
$$
If the moduli space $\underline 
{\hat \M}_T(S_{\a}, S_{\b})$ is empty, we define 
$\pa_{\a, \b}$ 
to be the zero operator. If it is nonempty, we define 
$\pa_{\a, \b}$ as follows. For $\sigma 
\not \in C_*(S_{\a})$, we set $\pa_{\a, 
\b} \sigma = 0$. For $\sigma=(\D_j,f)\in C_j(S_{\a})$ with 
$\mu(\a) + j  = k$, 
consider the fiber product 
$$\Delta=
 \D_j\times_{S_{\a}}\underline {\M}_T(S_{\a},
 S_{\b}) =\{(z,u)\in \D_j\times
 \underline {\M}_T(S_{\a},S_{\b}):f(z)=\pi_-(u)\}.$$

We have a natural map $\ti  \pi_+: \Delta \to 
S_{\b}, \ti \pi_+((z, u))=\pi_+(u).$
The fiber product $\Delta$ is a $(j+\mu(\a)-\mu(\b)-1)$-dimensional 
compact oriented manifold with corners, hence can be 
triangulated into oriented simplices which are identified with Euclidean
simplices. (We choose such a triangulation and identifications arbitrarily.)
The map $\ti \pi_+: \Delta \to 
 S_{\b}$ then gives rise to a singular $(j+\mu(\a)
 -\mu(\b)-1)$-chain in $S_{\b}$. We define this chain 
 to be $\pa_{\a, \b}\sigma$. 
 Clearly, we indeed have
 $$
 \pa_{\a, \b}: C_k \to C_{k-1}
 $$
for all $k$. 

\definition{Definition 7.4} We define $ \pa : C_k
\to C_{k-1}$ as follows. First, we define $\pa_0: C_k \to C_{k-1})$ by 
$ \pa_0= (-1)^k \pa_{O}$. Then we set $\pa = \pa_0  + 
\sum_{\mu(\a) > \mu(\b)} \pa_{\a, \b}$. 
Next, we define  $\ti \pa: \ti C_k \to \ti C_{k-1}$
by $\ti \pa <\sigma> = <\pa \sigma>.$ The 
boundary operator $\ti \pa: \ti C_*(\R^0) \to
\ti C_*(\R^0)$ is defined to be the direct sum of these boundary operators.
  (We abuse notations a bit here in order to avoid 
too many notations.) 
\enddefinition

\proclaim{Lemma 7.5} $\ti \pa^2=0$. 
Hence $(\ti C_*, \ti \pa)$ is a chain complex.
 \endproclaim

 \demo{Proof}
Consider $\sigma = (\Delta_j, f)\in C_j(S_{\a})$ 
 We have $\ti \pa^2 \sigma = \sum_{\mu(\b)  < \mu(\a)} I_{\b}$, where 
 $$
 I_{\b} = <\pa_0 \pa_{\a, \b}\sigma>
 + <\pa_{\a, \b}\pa_0\sigma>
 + \sum_{\mu(\a) > \mu(\ga) > 
 \mu(\b)}<\pa_{\ga, \b} \pa_{\a,
 \ga} \sigma>.  
$$
On the other hand, by (7.1) we have
$$\eqalign{ \pa(\D_j  \!\!\times_{S_{\a}} \!\underline 
{\hat \M}_T(S_{\a},S_{\b})) &=
\pa \D_j\!\!\times_{S_{\a}}\underline {\hat \M}_T(S_{\a},S_{\b})\!+ \!
(-1)^{j+dim S_{\a}}\!\D_j \!\times_{S_{\a}}\!\!\pa \underline 
{\hat \M}_T(S_{\a},S_{\b}) \cr
&= 
\pa \D_j\times_{S_{\a}}\underline {\hat \M}_{T}(S_{\a},S_{\b})
+\cr
& (-1)^{j+ \mu({\a})}\!\!\!\!\!\!\!\!\!\sum_{\mu(\a) 
>\mu(\ga) > \mu(\b)} \!\!\!\!\!\!\!\D_j \times_{S_{\a}} \underline 
{\hat \M}_T(S_{\a},S_{\ga}) \times_{S_{\ga}}
\underline {\hat \M}_T(S_{\ga},S_{\b}).\cr} \tag 7.4 
$$
Multiplying this equation by $(-1)^{(j+ \mu(\a) -\mu(\b) -1) + \mu(\b)}$ 
and noting Remark 7.3 we then infer  that $I_{\b}=0$. The desired 
identity follows. (This is similar to the situation in $\FU$.)
\qed 
\enddemo

\definition{Definition 7.6} We define the Bott-type Seiberg-Witten 
homology $FH^{SW}_{b*}(c)$ and cohomology $FH^{SW*}_{b}(c)$ to be the 
homology $H_*(\ti C_*, \ti \pa)$ and cohomology $H^*(\ti C_*$ $,
\ti \pa)$ of the complex $(\ti C_*, \ti \pa)$.
 (Recall that 
$c$ is the given $spin^c$ structure.)
\enddefinition 

\definition{Remark 7.7} In the above 
situation and in Section 9 of Part II, we can enlarge 
the space of singular chains $C_*(\X)$ to 
allow all $(\Delta, f)$, where $\Delta$ is a compact oriented 
smooth manifold  with corners and $f$ a continuous map from $\Delta$ 
to $\X$. Then we can avoid baricentric subdivisions and 
triangulations (we still divide out homotopy differences) 
and obtain equivalent homology and cohomology theories. 
\enddefinition

\head  8. Invariance I
\endhead

Consider two chain complexes $(C_*, \pa)$ and 
$(\bar C_*, \bar \pa)$. A chain map of degree $m \in \ZZ$  between 
the two complexes consists of homomorphisms $F: C_k \to \bar C_{k+m}$ 
such that $F \cdot \pa = \bar \pa \cdot F$. 
A shifting homomorphism $F: H_*(C_*, \pa) 
\to H_*(\bar C_*, \bar \pa)$ of degree $m \in \ZZ$ consists 
of homomorphisms $F: H_k(C_*,\pa) \to H_{k+m}(\bar C_*, \bar \pa)$. 
One defines shifting homomorphisms for cohomologies 
in a similar way.

Our goal is to prove the following invariance result.

\proclaim{Main Theorem I} 
The Bott-type Seiberg-Witten Floer 
homology and cohomology are diffeomorphism 
invariants  modulo shifting isomorphisms. 
\endproclaim 

In othere words, these homology and cohomology are independent of 
the metric $h$ and generic parameter $(\pi, \lmd)$
modulo shifting isomorphisms.

In this section, we construct the first kind of shifting homomorphisms 
which we need. In Part II, we present a second kind of shifting 
homomorphisms. We shall show that they provide  inverses for each other.

Consider two metrics $h_+$ and $h_-$ on $Y$ and 
generic parameters $(\pi_+,\lmd_+)$ for $h_+$ and 
$(\pi_-, \lmd_-)$ for $h_-$ respectively. We proceed to 
construct shifting homomorphisms between  our 
homologies (cohomologies) constructed with respect to $(h_+, \pi_+,
\lmd_+)$ and $(h_-, \pi_-, \lmd_-)$
respectively. 

Choose a smooth path of metrics $h(t)$ on  $Y$  such that
$$h(t)=\cases h_-, & \hbox{if }t<-1,\cr
h_+, & \hbox{if }t>1,\cr\endcases
$$
a smooth path of $\pi(t)\in\Pi$ 
$$\pi(t)=\cases \pi_-, & \hbox{if }t<-1,\cr
\pi_+, & \hbox{if }t>1,\cr\endcases
$$
and a smooth function $\lmd(t)\in \RR$ such that
$$\lmd=\cases \lmd_-, & \hbox{if }t<-1,\cr
\lmd_+, & \hbox{if }t>1.\cr\endcases
$$

The following lemma is an immediate consequence of 
Lemma C.1 in Appendix C.

\proclaim{Lemma 8.1} Fix an $(A_0, \P_0)
\in \A(X) \times \Ga^+(X)$. For any $R>0$, set $X_R= Y \times 
[-R, R]$ and
$$\eqalign{\SS_R=\{&u = (A_0, \P_0)+ (A, \P): (A, \P) \in \O^1_1
(X_R) \times \Ga^+(X_R)\cr
&\hbox{with }d^* A=0\hbox{ and }A|_{\pa X_R}(\dt) =0. \}\cr}$$
Then $\SS_R$ is a global slice for the action of $\G_2^0(X_R)$
on $\A_1(X_R) \times \Ga_1^+(X_R)$. In other words, $\A_1(X_R)
\times \Ga_1^+(X_R) = \G^0_2(X_R)\cdot \SS_R$.
\endproclaim

The following definition is a crucial construction. 

\definition{Definition 8.2} Choose a nonzero $\Psi_0 \in \Ga^-(X)$ with 
support contained in the interior of $X_R$. We define 
a smooth vector field $Z$ on $\A_{1}(X_1) \times \Ga^+_{1}(X_1) $ by
$$
Z(g^*u )=g^{-1}\Psi_0,
$$
for $g \in \G^0_2(X_R), u \in \SS_R$ and extend $Z$  to $\A_{1, loc}(X) 
\times \Ga^+_{1, loc}(X)$ by 
$$
Z(u)=Z(u|_{X_1}).
$$
\enddefinition
We endow $X$ with 
the warped product metric 
determined by the family of metrics $h(t)$ and the standard metric 
on $\RR$.

The following lemma is readily proved. 
\proclaim{Lemma 8.3}
$Z$ is equivariant with respect to the action of $\G_{2, loc}^0$. 
\endproclaim

Now we introduce the (perturbed) transition trajectory 
equation  for $A= a + f dt, \P=\p$. 
$$
\eqalign{
\frac{\partial a}{\partial t} & =*F_a+d_Yf+\langle e_i\cdot\p,\p\rangle e^i+
\n H_{\pi(t)}(a) + \e b_0,
\cr
\frac{\partial \p}{\partial t} & =-\np_a\p-\lmd(t)\p+\e Z,\cr
}
\tag 8.1
$$
with additional parameters $ \e \in \RR$ and $b_0$, which is 
a smooth 1-form of compact support (it does not contain $dt$). 
The Hodge $*$ at time $t$ in the 
equation is that of the metric $h(t)$. 
The perturbation term $\e  Z$ is called a {\it spinor 
perturbation}.  We have  the following obvious, but crucial lemma. 

\proclaim{Lemma 8.4} The equation (8.1) 
is equivariant  with respect to the action of $\G^0_{2, loc}$.
Moreover, it has no reducible solution.
\endproclaim

A fundamental property of the Seiberg-Witten equation 
is a pointwise maximum principle  for the spinor 
field, which is a consequence of 
the Weitzenb{\"o}ck formula (2.1). With the presence of $Z$, this 
principle no longer holds. Instead, we have the following 
result.

\proclaim{Lemma 8.5} Let $(A, \P)=(a+ f dt, \p) $ be 
a solution of (8.1). Then there holds
$$
\|\P\|_{L^{\infty}} \leq C E(A, \P)
$$
for a constant $C$ depending only on the 
families $h(t), \pi(t), \lmd(t)$ and the geometry of $Y$. 
\endproclaim
\demo{Proof} Before proceding 
with the proof, we first observe that 
by (2.11) and (8.1) the energy can be estimated in the following way
$$
\eqalign {E(A,\P)=2 \lim_{t \to
-\infty} \hbox{ \bf cs}_{(\lambda, H)}(a(\cdot, t), \p(\cdot, t))
\cr -2\lim_{t \to +\infty} \hbox{\bf cs}_{(\lambda, H)}(\a(\cdot,
t), \p(\cdot, t))+\int_X|Z|^2, \cr} \tag 8.2$$ 
$$\int_X|Z|^2<C. \tag 8.3 $$
Using local Columb gauges provided 
by Lemma C.1 in Appendix C and a patching argument,
we can perform a gauge transformation to convert $(A, \P)$ into a smooth 
solution. Since the $L^{\infty}$ norm of $\P$ is invariant, 
we can assume that $(A, \P)$ is already smooth. Furthermore, 
we can assume that $(A, \P)$ is in temporal form.

For simplicity, we assume $\lmd(t) \equiv 0$ in the following 
argument. It is easy to modify it to handle $\lmd(t)$. Put
$$I_1= \np_a\p, I_2= *F_a+\langle e_i\cdot\p,\p\rangle
e^i - \n  H(a, \p), I_3=I_2 + \n H(a, \p).$$
Then 
$$\int_{\Omega}(|I_1|^2+|I_2|^2) = \frac 12 E(A, \P, \Omega),$$
for any domain $\Omega \subset X$.
For each $t\in\RR$, we use the 3-dimensional  Weitzenb\"ock formula (2.1) 
on $Y\times\{t\}$ to derive 
$$\np_aI_1=-\D_a\p+\frac{\bar s}4\p+\frac 12|\p|^2\p+I_3\p,\tag 8.4$$
where $\bar s$ denotes the scalar curvature function of $(Y, h(t))$.
Multiplying (8.4) by $\p$ and integrating by parts, we infer  
$$\int_{Y\times\{t\}}(|\n_a\p|^2+\frac s4|\p|^2+ \frac 12 |\p|^4) \leq  
\int_{Y\times\{t\}}|I_1|\cdot|\pa_a\p|+|I_3||\p|^2.$$
Using the H\"older ineqality we then  deduce that
$$\int_{Y\times\{t\}}(|\n_a\p|^2+ \frac 14
|\p|^4) \leq C(1+\int_{Y\times\{t\}}(|I_1|^2+|I_2|^2),$$
where $C$ depends only on $\|\n H\|_{L^{\infty}}$, which 
can be estimated by appealing to Lemma 3.11.
This last estimate implies 
$$\int_{X_{R-2, R+2}}(|\n_a \p|^2+ |\p|^4) \leq C(1 + E(A, \P, 
X_{R-2, R+2}))
\tag 8.5$$
for any $R >0$, where $X_{r, R} =Y \times [r, R]$. 

Next we apply the Moser iteration to deduce 
the desired $L^{\infty}$ estimate. Let  
$\xi:X\to [0,1]$ be a cut-off function such that
${supp~}\xi \subset X_{R-2, R+2}$ and $\xi(t)=1$ 
for $t\in X_{R-1, R+1}$. 
By the 4-dimensional  Weitzenb\"ock formula (2.1) on $X$ (recall 
that $X$ is endowed with warped product metric), we have
$$
D_AZ=-\D_A \P+\frac{s}4 \P-\frac14 |\P|^2\P.
\tag 8.6$$
Choosing $\xi^2|\P|^p\P$ as a test function,
where $p>0$ will be determined later, we obtain 
$$\eqalign{(\int_X\xi^2\n_A\P\cdot\n_A(|\P|^p\P)+
\frac{s}4|\P|^{p+2}+\frac14 |\P|^{p+4})\cr
=-\int_X(\xi^2 ZD_A(|\P|^p\P)+2\xi\n \xi Z|\P|^p\P+
2\xi\n \xi\cdot\n \P|\P|^p\P.\cr}
\tag 8.7$$

We have 
$$
\eqalign{\n_A(|\P|^p\P) &=|\P|^p\n_A\P+d|\P|^p\P \cr
&=|\P|^p\n_A \P+p|\P|^{p-1}\frac{\langle \n_A\P,\P\rangle}{|\P|}\P,
\n_A\P\cdot\n_A(|\P|^p\P) \cr &=|\P|^p|\n_A\P|^2+p|\P|^{p-2}\langle\n_A\P,
\P\rangle^2. \cr}$$
On  the other hand,
$$ p \|\P\|^{p-2}  \langle \n_A \P, \P  \rangle^2 \leq C \varepsilon
\|\P\|^p \|\n_A \P\|^2 + \frac {Cp^2}{\varepsilon}  \|\P\|^p,
$$
where $\varepsilon > 0$ is arbitrary. Similarly,
$$\eqalign{|D_A(|\P|^p\P)| & \le|\P|^p|D_A\P|+|d|\P|^p\cdot\P| \cr
& \le C |\P|^p|\n_A\P|+p|\P|^p|\n_A\P| \cr
& \le C \varepsilon  |\P|^p|\n_A\P|^2+\frac C{\varepsilon} (p^2+1) |\P|^p,
\cr}$$

$$|\n \xi||\xi||Z||\P|^{p+1}\le C\xi^2|\P|^{p}
+C|\n\xi|^2|\P|^{p+2},$$
and
$$|\n \xi||\xi||\n
\P||\P|^{p+1}\le\varepsilon\xi^2\|\P\|^p |\n\P|^2+ \frac C{\varepsilon}
|\n \xi|^2|\P|^{p+2}.$$
Choosing $\varepsilon$ suitably, we deduce 
$$\int\xi^2|\P|^p|\n_A\P|^2\le C\int((p^2 +1)\xi^2|\P|^p+|\n
\xi|^2|\P|^{p+2}).$$
Consequently,
$$\int_X \xi^2|\P|^p|\n|\P||^2\le C(\int_X(p^2+1) \xi^2|\P|^p+
|\n \xi|^2|\P|^{p+2}),$$
or
$$
\int_X\xi^2|\n|\P|^{\frac{p+2}2}|^2\le C(p+1)^2\int_X(\xi^2|\P|^p+|\n
\xi|^2|\P|^{p+2}).
$$
Now we set $w=|\P|^{(p+2)\slash 2}$.
By the Sobolev inequality and  H\"older inequality, we arrive at 
$$\eqalign{\|\xi w\|_{L^2}^2 & \le C(p+1)^2\int 
(|\xi\n w|^2+|w\n \xi|^2)\cr
&\le C(p+1)^2(\|\xi
w\|_{L^2}^{\frac{p}{p+2}}+\||\n\xi|w\|_{L^2}).\cr}$$
Then we use the iteration process as presented in $\Bg$ to infer 
$$\sup_{X_{R-1, R+1}}|\P|\le C\|\P\|_{L^4(X_{R-2, R+2})}.$$
Combining it with (8.5) we are done. 
\qed
\enddemo

We also have 
\proclaim{Proposition 8.6} An analogue of Proposition 6.2 for (8.1) holds.
\endproclaim

We have various configuration spaces  and moduli spaces 
associated with (8.1) which are analogous to the spaces introduced in 
Section 4. All the analysis in Sections 4, 5 and 6 
carres  over.  We shall be brief in formulating the relevant results.

Let e.g. $\R_{\pm}$ denote the $\R$ for $(h_{\pm}, \pi_{\pm}, \lmd_{\pm})$.
Consider 
$\a_- \in \R_{-} $, $\a_+ \in\R_{+} $ and $p_{\pm}\in\a_{\pm}$. 
We have the space of transition trajectories $\N(p_-, p_+)$ 
and the moduli spaces $\M^0(p_-, p_+), \M_T(S_{\a_-},  
S_{\a_+})$ etc..
We also have the various spaces of transition $k$-trajectories.
A  transition $k$-trajectory is a $k$-tuple $(u_1,...,u_k)$  with 
a distinguished portion $u_m, 1 \leq m \leq k$, such that 
$u_i \in \N(p_{i-1}, p_i)$ with $p_0=p_, p_k=p_+;
p_i \in \SS\W^{-}, 0 \leq i \leq m-1; p_i 
\in \SS\W^{+}, m \leq i \leq k.$  Now ``proper" means 
that $p_0,...,p_{m-1}$ belong to distinct 
gauge classes, and $p_{m},...,p_k$ belong to distinct 
gauge classes. The twisted time translation acts on all portions except 
the distinguished one. The other concepts regarding 
$k$-trajectories carry over easily to transition $k$-trajectories. 

We have the following analogues of Proposition 4.12, 
Remark 4.14 and Proposition 6.13. 

\proclaim{Proposition 8.7} Let $\e$ be given. Then for generic 
$b_0$ $($we shall say that $(\e, b_0)$ is generic$)$, transversality holds for 
the moduli spaces $\M^0(p_+, p_-)$ with $p_{\pm}
\in \a_{\pm} \in \R_{\pm}$. Consequently, it 
is a smooth manifold of dimension $ \hbox{ind }\F_{p_-, p_+} -
\hbox{max }\{\hbox{dim }$ 
$G_p, \hbox{dim }G_q\} +1$, where $ \F_{p_-, p_+}$ 
means the linerization of the operator in (8.1) (the 
{\it transition Seiberg-Witten operator}). Moreover,  
the spaces
$\underline {\M^0}(p_0,...,p_k)$ with $ p_0 \in \SS\W^{-}_0,
p_k \in \SS\W_0^+$ are smooth manifolds. 
\endproclaim

\proclaim{Proposition 8.8} For generic $(\e, b_0)$, we have 
for all $\a_- \in \R_{-}, \a_+ \in \R_{+}$

(1) ${\underline {\hat M}}_T(S_{\a_-},
S_{\a_+})$ has the structure of
$d$-dimensional smooth oriented manifolds with
corners, where  
$d = \hbox{ind }\F_{p_-, p_+} -\hbox{dim }G_p- \hbox{dim }G_q +1$
for $p_- \in \a_-, p_+ \in \a_+$.

(2) This structure is compatible with the natural stratification of 
${\underline {\hat \M}}_T(S_{\a_-},
S_{\a_+})$. 

(3) The temporal endpoint maps $\pi_-:
{\underline {\hat \M}}_T(S_{\a_-},
S_{\a_+}) \to  S_{\a_-}$ and $\pi_+: {\underline {\hat \M}}_T(S_{\a_-},$
$S_{\a_+})$ $ \to S_{\a_+}$ are smooth maps. But they 
are not fibrations in general. 
\endproclaim

The proof of the first statement of Proposition 8.6 is analogous to 
the proof of Proposition 4.13 and Remark 4.15, because Lemma 8.4 rules out 
reducible transition trajectories. Note that 
instead of using holonomy perturbations we now use 
the perturbation $b_0$ as in $\KM$. (This perturbation 
is not time translation equivariant, and hence can't 
be applied in the construction of our homologies.)

Now we proceed to construct our first kind of shifting homomorphisms. Let 
$O_{\pm}$ be the unique reducible elements in $\R_{h_{\pm}}$ respectively. We set 
$$
m_0= {ind~ }\F_{O, O'} -1. \tag 8.8
$$

The following lemma is analogous to Corollary 5.2.
\proclaim{Lemma 8.9} We have 
$$
{ind }~\F_{p_-, p_+} = \mu_-([p_-])-\mu_+([p_+])
+m_0 + {dim }~ G_{p_+}.
\tag 8.9
$$
Consequently,
$$
{dim }~ \M_T(S_{\a_-}, S_{\a_+}) = \mu_-(\a_-)
-\mu_+(\a_+)
+m_0 - {dim } ~G_{\a_-}+1.
\tag 8.10
$$
\endproclaim

Now we define  homomorphisms $F: \ti C_{k}(\R_-) \to 
\ti C_{k+m_0}(\R_+). $ Consider $<\sigma> \in C_j(S_{\a_-})$
with $\sigma = (\Delta_j, f)$. 
We choose a representative $\sigma$ such that 
$f$ and $f|_{\pa \D_j}$ are transversal 
to the endpoint maps $\pi_-$ from the  moduli spaces $\underline 
{\hat \M}_T(\a_-, \a_+)$ and their boundaries 
for all   $ \a_+ \in \R_+$.
For each $\a_+ $ with the corresponding 
moduli space nonempty, we follow the construction of the 
boundary operator $\pa_{\a,\b}$ in Section 7
to obtain a singular chain $\sigma' 
\in C_{j'}(S_{\a_+})$ with $j'+ \mu_+(\a_+) =
j+ \mu_-(\a_-) +m_0$. 
We define $F_{\a_-, \a_+}(<\sigma>)$ to 
be $<\sigma'>$. We define it to be zero 
if the moduli space is empty. Then we set 
$$
F(<\sigma>) = \sum_{\mu_+(\a_+) \leq \mu_-(\a_-)+m_0}
F_{\a_-, \a_+}(\sigma).
$$

It is easy to see that  we indeed 
have $F: \ti C_{k}(\R_-) \to \ti C_{k+m_0}(\R_+)$. 

\proclaim{Proposition 8.10} We have 
$\ti \pa \cdot F = F \cdot \ti \pa$, hence $F$ is a  chain map
of degree $m_0$ from $ (\ti C_*(\R_-), \ti \pa)$
to $(\ti C_*(\R_+), \ti \pa)$. The induced shifting 
homomorphisms between the homologies and cohomologies are 
denoted by $F_*$ and $F^*$.
\endproclaim
\demo{Proof} 
The proof goes along the same lines as the proof of Lemma 7.3. 
To simplify notations, 
we argue here in terms of $\underline {\hat \M}_T(S^-_i, S_k^+)
= \cup_{\a_- \in S^-_i, \a_+ \in
S^+_k} \underline {\hat \M}_T(\a_-, \a_+)$. 
Analogous to (7.4) (the moduli spaces $\underline {\hat \M}_T(
\a_-, \a_+)$ are oriented in a way similar to Lemma 7.2) we have 
for $\sigma=(\Delta_j, f)  \in C_j(S_i^-)$ and $k \leq i+m_0$
$$\eqalign{ \pa(\D_j  \times_{S_{i}^-} \underline
{\hat \M}_T(S_{i}^-,S_{k}^+)) &=
\pa
\D_j\times_{S_{i}^-}\underline {\hat \M}_T(S_{i}^-,S_{k}^+)+ \cr
(-1)^{j+dim S_{i}^-}\D_j \times_{S_{i}^-}\pa \underline
{\hat \M}_T(S_{i}^-, S_{k}^+)
&=
\pa \D_j\times_{S_{i}^-}\underline {\hat \M}_{T}(S_{i}^-,S_{k}^+)
+\cr}$$ $$ (-1)^{j+ i+ m_0 +1}\sum_{i
>m \geq  k-m_0} \D_j \times_{S_{i}^-} \underline
{\hat \M}_T(S_{i}^-,S_{m}^-) \times_{S_{m}^-}
\underline {\hat \M}_T(S_{m}^-,S_{k}^+)+  $$
$$ (-1)^{j+ i+ m_0+1}\sum_{i
\geq m'-m_0  > k -m_0} \D_j \times_{S_{i}^-}
\underline
{\hat \M}_T(S_{i}^-,S_{m'}^+) \times_{S_{m'}^+}
\underline {\hat \M}_T(S_{m'}^+,S_{k}^+).$$

Clearly, this implies
$$ <\pa_0 F_{S_i^-, S^+_k} \sigma > =
-<F_{S_i^-, S^+_k} \pa_0 \sigma> - \sum_{i>m \geq k-m_0}
<F_{S_m^-, S_k^+} \pa_{S_i^-, S^-_m} \sigma > -$$
$$
<\pa_{S^+_{m'}, S^+_k} F_{S_i^-, S^+_{m'}} \sigma>,
$$
where e.g. $\pa_{S_i^-, S^-_m}$ and $F_{S_m^-, S_k^+}$
are defined analogously to $\pa_{\a, \b}$ and 
$F_{\a_-, \a_+}$.  Summing over all $k$, we arrive 
at the desired chain homotopy property.
\qed
\enddemo

\head{9. Stable Bott-type and stable equivariant homology}
\endhead

\noindent {\bf Stable Bott-type Seiberg-Witten 
Floer homology} 

\bigskip

We consider
    $$C^s_k=\oplus_{i+j=k}C_j(S_i \times S^1),
     C^k_s=\oplus_{i+j=k}C^j(S_i)\times S^1,$$
     and the corresponding $\ti C^s_k, \ti C^k_s$,
     which are analogous to $\ti C_k, C^k$.

      Now the moduli spaces $\underline {\hat M}_T
       (S_{\a}, S_{\b})$ are replaced by
	$\underline {\hat M}_T(S_{\a}, S_{\b}) \times
	 S^1$, with endpoint maps:
	  $$
	   \pi_+(u, s)= (\pi_+(u), s),  \pi_-(u, s)= (\pi_-(u),s).
	    $$
	     We have  the boundary operator $\ti \partial_{s}:
	      \ti C^s_k \to \ti C^s_{k-1}$ 
	      analogous to $\ti \pa$.(Note that 
the moduli spaces of trajectories and the submanifolds 
$S_{\a} \times S^1$ are both one dimension higher 
than in Section 7. The extra dimensions 
cancel each other in the construction of the boundary operator.)
Similarly, we have $\ti \partial_{s}^2 =0$.
Hence we can introduce the following stable Seiberg-Witten homology
and cohomology of Bott type.

\definition{Definition 9.1} We define $FH^{SW}_{sb*}(c)
 =H_*(\ti C^s_*, \ti \pa_s)$ and $
FH^{SW*}_{sb}(c) = H^*(\ti C^s_*, \ti \pa_s)$. 
\enddefinition

 \bigskip

 \noindent {\bf Stable equivariant Seiberg-Witten
		       Floer homology}
		       \bigskip

			There is a diagonal action of $S^1$
			 on $\R^0 \times S^1$:
			  $$
			   s^*(\a, s')= (s^*\a, s^{-1}s'),$$
			    where $s ,s' \in S^1$, $s$ is 
			    identified with the constant map 
			    from $Y$ into $S^1$ with 
			    value $s$ (a constant gauge),
			    and $s^*[p]_0=[s^*p]_0$.
This action induces an action of $S^1$ on
singular chains in $\R^0 \times S^1$:
 if $\sigma=(\Delta, f)$ is a $j$-chain, and $s \in S^1$,
  then $(s^*\sigma )(z)=(\Delta_j, s^*f(\cdot))$.
   Passing to quotients, we
    obtain an action of $S^1$ on
     $\ti C_*(\R^0 \times S^1)$.  A
      $j$-cochain class
       $\omega  \in \ti C^*(\R^0
	\times S^1)$ is called {\it equivariant}
	 (or {\it invariant}),
	  provided that $\omega(s^*<\sigma>)=\omega(<\sigma>)$ for
	  all $<\sigma> \in \ti C_*(\R^0)$ and $s \in S^1$.

	  We define $\ti C^j_e(S_i \times S^1)$
	  to be the free abelian group of
	  equivariant $j$-cochain classes on
	  $S_i \times S^1$ and set

	  $$\ti C^k_e=\oplus_{i+j=k}C^j_e(S_i \times S^1).$$

	  Now let $\ti \partial^*_e$ be the restriction of
	  $\ti \partial^*_{s}$ to equivariant cochain
classes. 

\proclaim{Lemma 9.2} The endpoint maps $\pi_{\pm}$ 
defined on the moduli spaces $\underline {\hat \M}_T
S_{\a}, S_{\b}) \times S^1$ 
are $S^1$-equivariant, where the $S^1$ action on  
$\underline {\hat \M}_T(S_{\a}, S_{\b}) \times S^1$ is induced from 
the following action: $s \in S^1$ acts on $(u, s')$ to 
yield $(s^*u, s^{-1}s')$.  
(If we work with the model $\underline {\hat \M^0}(p, q)$, 
then the $S^1$ action is induced from that on 
$\underline {\hat \M}_T(S_{\a}, S_{\b})$ through 
the temporal transformation.)
\endproclaim 

We omit the easy proof. 
As a consequence of this lemma, 
the operator $\ti \pa_e^*$ has equivariant 
cochain classes as values. 
Hence we have 
$\ti \partial^*: \ti C^{k}_e \to \ti C^{k+1}_e$.

\definition{Definition 9.3} The stable
equivariant Seiberg-Witten Floer
cohomology
$FH^{SW*}_{se}$
$(c)$ is defined to be the 
homology $H_*(\ti C^*_e,
\ti \pa^*_e)$ of
the complex $(\ti C^*_e,
\ti \pa^*_e)$. 
The stable equivariant Seiberg-Witten
Floer homology
$FH^{SW}_{se*}(c)$ is defined
to be the cohomology
$H^*(\ti C^*_e, \ti \pa_e^*)$
of the complex $(\ti C^*_e, \ti \pa_e^*)$.
\enddefinition

\head{Invariance II}
\endhead

The purpose of this section is to prove 
the following two results.

\proclaim{Main Theorem II} The stable 
Bott-type Seiberg-Witten Floer homology
and cohomology are diffeomorphism 
invariants up to shifting isomorphisms.
\endproclaim

\proclaim{Main Theorem III} The stable equivariant Seiberg-Witten 
Floer homology and cohomology are diffeomorphism invariants 
up to shifting isomorphisms.
\endproclaim

We start with 

\definition{Definition 10.1}
Let $\G_{2, loc}$ act on $(\A_{1, loc}
\times \Ga^+_{1, loc}) \times S^1$
in the following fashion:
$$
(gg_0)^*(u, s) = ((gg_0)^*u,
g_0^{-1}s),
$$
where $g \in S^1$ (the group of constant gauges), $g_0\in
\G_{2, loc}^0, $ and $ s \in S^1$.
\enddefinition

\definition{Definition 10.2} We define a smooth vector field
$Z_e$ on $(\A_{1, loc} \times\Ga^+_{1, loc}) \times S^1$
as follows $$Z_e(u, s) = s Z(u), $$
where $Z$ is the vector field given by Definition 8.3.
\enddefinition

The following lemma is readily
proved.
\proclaim{Lemma 10.3} $Z_e$ is equivariant with respect to
the action of $\G_{2, loc}$.
\endproclaim

Now we introduce the following stable version of 
the transition trajectory equation  (8.1)  for $A= a + f dt, \P=\p
$ and $ s \in S^1$. Its solutions will be called {\it stable transition 
trajectories}.

$$
\eqalign{
\frac{\partial a}{\partial t} & =*F_a+d_Yf+\langle e_i\cdot\p,\p\rangle e^i+
\n H_{\pi(t)}(a) + \e b_0,
\cr
\frac{\partial \p}{\partial t} & =-\np_a\p-\lmd(t)\p+\e sZ_{e}(A, \P),\cr
}
\tag 10.1
$$
where $s \in S^1$ and $\e, b_0$ are the same as in (8.1).

Now consider $(h_{\pm}, \pi_{\pm}, 
\lmd_{\pm})$ as in Section 8. 
For $\a_{\pm} \in \R_{\pm}$ and $p_{\pm} \in \a_{\pm}$, 
we consider $\a_{\pm} \times  S^1$ and set 
$$
\N(p_-\times S_1, p_+ \times S^1) =\{ (u, s): 
(u, s) \hbox{ is a stable transition 
trajectory converging to}
$$
$$ p_{\pm} \hbox{ at } 
\pm \infty. \}
$$

We have moduli spaces $\M^0(p_-\times S^1, p_+ \times 
S^1), \M_T(S_{\a_-}\times S^1,
S_{\a_+} \times S^1)$ etc.. 
We also have the various spaces of stable transition $k$-trajectories.
A  stable transition $k$-trajectory is a pair $((u_1,...,u_k),s)$  with
$s \in S^1$ and
a distinguished portion $u_m, 1 \leq m \leq k$, such that
$(u_i, s)  \in \N(p_{i-1}\times S^1, p_i \times S^1)$ with
 $p_0=p_-$,$ p_k=p_+;
p_i \in \SS\W^{-}, 0 \leq i \leq m-1; p_i
\in \SS\W^{+}, m \leq i \leq k.$  Now  ``proper" means
that $p_0,...,p_{m-1}$ belong to distinct
gauge classes, and $p_{m},...,p_k$ belong to distinct
gauge classes. The twisted time translation is easily defined and 
acts on all portions except
the distinguished one. The concept of consistent 
stable transition $k$-trajectories is also easily defined. 
Finally, note that $\pi_{\pm}(u, s)$ is defined to be $(\pi_{\pm}(u), s)$.

We have the following analogues of Proposition 4.12,
remark 4.14 and Proposition 6.13.

\proclaim{Proposition 10.4} Let $\e$ be given. Then for generic
$b_0$ (we shall say that $(\e, b_0)$ is generic), transversality holds for
the moduli spaces $\M^0(p_+ \times S^1, p_- \times S^1)$ 
for $p_{\pm} \in \a_{\pm} \in \R_{\pm}$. Consequently, it
is a smooth manifold of dimension $ \mu_-(\a_-) - 
\mu_+(\a_+) +m - \hbox{ dim } G_{\a_-} + 1$. 
 Moreover,
the spaces
$\underline {\M^0}(p_0 \times S^1, ..., 
p_k \times S^1)$ with $p_0 \in \a_-, 
p_k \in \a_+$ are smooth manifolds. 
\endproclaim

\proclaim{Proposition 10.5} For generic $(\e, b_0)$, we have
for all $\a_- \in \R_{h_-}, \a_+ \in \R_{h_+}$

(1) ${\underline {\hat M}}_T(S_{\a_-} \times S^1,
S_{\a_+}\times S^1)$ has the structure of
$d$-dimensional smooth oriented manifolds with
corners,
where $d = \mu_{-}(\a_-)- \mu_+(\a_+) +m -\hbox{ dim }
G_{\a_-} +1.$

(2) This structure
is compatible with the natural stratification of
$\underline {\hat \M}_T(S_{\a_-} \times S^1,
S_{\a_+}\times S^1)$. 

(3) The temporal endpoint maps $\pi_{\pm}:
{\underline {\hat \M}}_T(S_{\a_-} \times S^1,
S_{\a_+}\times S^1) \to  S_{\a_{\pm}} \times S^1$
are $S^1$-equivariant smooth maps, where the action of $S^1$ on the 
moduli spaces is induced by the following action: $s \in S^1$ 
acts on $(u, s')$ to yield $(s^*u, s^{-1}s')$.
But they are not fibrations in general.
\endproclaim

To prove Main Theorem II, we construct chain maps $F^{+}:\ti C^{s+}_* \to 
C^{s-}_*$ as in Section 8. We need to show that they induce isomorphisms 
between the homologies and cohomologies. The arguments consist of three  
steps. 
\bigskip

\noindent {\bf Step 1} 

Consider the equation (10.1) with parameters 
determined by the construction of $F^{-}$ and the same equation with 
parameters determined by the construction of $F^+$. We interpolate 
(glue) the former with the latter to obtain a new equation of a similar 
type, which we call the {\it glued transition equation}.
There is an interpolation parameter $\rho$ in this equation. 
When $\rho \to 0$, the equation breaks into the original two equations.
Fix some $\rho_0 >0$. Using the 
glued transition equation with this parameter we construct 
a chain map $F^-_+: \ti C^-_* \to \ti C^-_*$ of degree zero in the same 
way as e.g. the construction of $F^-$. (The compactified moduli spaces 
for this equation has the same structure as described in Proposition 10.7.)
Employing the enlarged glued transition equation in which the parameter 
$\rho$ varies in the range $0< \rho \leq \rho_0$, we then obtain 
a chain map $\Theta: \ti C^-_* \to\ti C^-_*$ of degree one 
(instead of zero, because
the corresponding moduli spaces are one dimension higher 
than before). (The compactified moduli spaces 
have a similar structure to that described in Proposition 10.7.)
Arguing as in the proofs of Lemma 7.2 and Proposition 8.10,
we infer that this map is a chain homotopy between 
$F^+ \cdot F^-$ and $F^-_+$, i.e. $F^+ \cdot
F^- - F^-_+ = \ti \pa \cdot \Theta + \Theta \cdot
\ti \pa. $

\bigskip

\noindent {\bf Step 2}

Next we make use of the parameter $\e$ in (8.1), 
which enters  the glued transition equation.  (There are actually two 
such parameters in the glued transition equation. We can turn  
them into one parameter.) When $\e =0$, the glued transition equation 
reduces to the trajectory equation (4.1) with some holonomy 
and $\lmd$ perturbations. We can incorporate additional 
holonomy perturbations to ensure transversality 
between irreducible Seiberg-Witten points 
or one irreducible and one reducible. The transversality 
between two reducibles follows from Lemma C.2.   Then we can 
construct a chain map $F_0^-: \ti C^-_* \to 
\ti C^-_*$ of degree zero. It is easy to see that 
it induces the identity isomorphism of the homology and cohomology.

Allowing $\e $ to vary from the given generic value to zero, we obtain  
another enlarged glued transition equation. Using it 
we then obtain a chain homotopy between $F^-_+$ and $F_0^-$.
It follows that $ F^+_* \cdot F^-_* = Id, F^{-*} \cdot F^{+*}=Id$.

\bigskip

\noindent {\bf Step 3} 

A  similar construction with the roles of $ F^-$ and $F^+$ 
reversed yields $F^-_* \cdot F^+_*=Id, F^{+*} \cdot F^{-*} =Id. $

We have proved Main Theorem II. The same arguments apply to 
Main Theorem III because of the equivariance of the constructions. 

\head Appendix A Seiberg-Witten Floer homology
\endhead

Here, we only consider rational homology spheres. 
It is not hard to extend the construction to general manifolds. 

Let $Y$ be a rational homology sphere 
with a given metric $h$ and $c$ a $spin^c$ structure on $Y$ as 
before. Set $C_i=\ZZ\{\a\in\R^*|\mu(\a)=i\}$.  We
define a boundary operator $\pa:C_i\to C_{i-1}$ in terms of the 
moduli spaces $\ti \M(\a, \b)$, or equivalently
$\ti \M(p, q), p \in \a, q\in \b$, 
where the tilde means quotient by the time translation action.
For a generic pair $(\pi, \lmd)$, $\ti \M(\a, 
\b)$ is a compact oriented manifold of zero dimension,
provided that $\mu(\a) - \mu(\b)=1$. The orientation (sign)
is given by pairing the orientation of $\M(\a, \b)$ 
provided by Proposition 5.3 with its orientation induced 
by the time translation action. For $\a \in C_i$ we then set 
$$\pa \a=\sum_{\mu(\a)-\mu(\b)=1}\sharp\ti\M(\a, \b)\b,\tag A.1$$
where $\sharp\ti\M(p,q)$ is the algebraic sum of $\ti\M(\a, \b)$.

The compactification of the moduli spaces $\ti \M(\a, \ga)$ with 
$\mu(\a) - \mu(\ga) =2$ is similar to 
the results in Section 6. For dimensional reasons, 
no trajectory connecting to 
the reducible point appears in the compactification.
Using these compactified moduli spaces and the 
consistency of orientation (Proposition 5.3) 
we obtain 

\proclaim{Lemma A.1} $\pa^2=0$.
\endproclaim

\definition{Definition A.2} The Seiberg-Witten Floer homology 
$FH^{SW}_*(c, h, \pi, \lmd) $ and 
cohomology $FH^{SW*}(c, h, \pi, \lmd) $ 
for the $spin^c$ structure $c$ and the parameters $h, \pi, 
\lmd$  are defined to be the homology and 
cohomology of the chain complex $(C_*, \pa)$. 
\enddefinition

\head{Appendix B Equivariant 
Seiberg-Witten Floer homology}
\endhead

There are two possible versions. The first 
uses equivariant singular cochains. The construction 
is similar to equivariant construction in Section 
9 of Part II. We use equivariant singular cochains 
on $S_{\a}$ instead of $S_{\a} \times 
S^1$. The resulting homology and 
cohomology will be denoted by $FH^{SW}_{e*}$ and 
$FH^{SW*}_e$.

The second version uses equivariant differential forms.
 Let $G$ be a compact Lie group,
 ${\frak g}$ its Lie algebra and $M$ a $G$-manifold.
 Let $\CC[{\frak g}^*]$ denote the algebra
 of complex valued polynomial function on ${\frak g}$.
 We  define the space of equivariant
  differential forms
  $$\O_G(X)=(\O^*(X)\otimes \CC[{\frak g}])^G$$
  to be the subalgebra of $G$-invariant elements in
  the algebra $\O^*(M)\otimes \CC[{\frak g}]$. The algebra $\O^*(M)\otimes
  \CC[{\frak g}])$
  has a $\ZZ$-grading defined by
  $$\deg(w\otimes z)=2\deg(z)+\deg(w),$$
  where $w\in \O(M)$ and $z\in \CC[{\frak g}]$.
We define   a differential $d_G $: 
  $$(d_G\a)(X)=d(\a(X))-\iota(X)(\a(X)),$$
  where $\a\in\O_G(M)$ and  
  $\iota(X)$ denotes the contraction by the vector field $X_M$ induced
  by an element $X\in{\frak g}$. We have $d_G^2=0$.

  To define the desired equivariant 
homology and cohomology for 
given generic parameters $(\pi, \lmd)$, we set 
  $$C^{i,j}=\O^j_{S^1}(S_i), C^k = \oplus_{i+j =k} C^{i, j},$$
  where $j$ denotes the degree of differential forms. (
The $S_i$ were defined in Section 7.)
  We define operators $\pa_r:C^{i,j}\to C^{i+r,j-r+1}$ by
  $$\pa_r\o=\cases d_{S^1}\o, & \hbox{if }r=0, \cr
  (-1)^j\pi_{+*}\pi_-^*(\o) & \hbox{otherwise},\cr
  \endcases $$
where the map $\pi_{+*}$ is integration along the fiber of the
 fibration $\pi_+: \underline {\hat \M}_T(S_i, S_{i+r})$ $ \to
S_{i+r}$ (cf. $\BA$), and $\pi_-$ refers to the fibration 
$\pi_-: \underline {\hat \M}_T(S_i, S_{i+r} ) \to S_{i+r}$.
($\underline {\hat \M}_T(S_i, S_j) \equiv \cup_{\a \in 
S_i, \b \in S_j} \underline {\hat \M}_Y(\a, \b)$. )  
 We set 
$$\pa_{S^1}=\sum \pa_r.$$

Using the compactification results in Section 6 
and the arguments in $\Ba$, one easily deduces 

  \proclaim{Lemma B.1} $\pa_{S^1}^2=0$.
  \endproclaim

  \definition{Definition B.2} The deRham type  
equivariant Seiberg-Witten Floer
  homology $FH^{SW}_{de*} (c, h, \pi, \lmd)$ 
for a given metric $h$ and generic parameters 
$h, \pi, \lmd$ is defined to be the 
homology of the chain complex $(C^*, \pa_{S^1})$.
The corresponding cohomology $FH^{SW*}_{de}$
is defined by a dual construction, cf. $\BA$. 
\enddefinition

\head{Appendix C Two analysis lemmas}
\endhead

First, we prove a result on local Columb gauge fixing. 
We assume that $Y$ is a rational homology 
sphere. Set $X_{r, R}=Y\times
[r ,R]$ for $r < R$ and fix a point $x_0 \in X_{r, R}$.  
\proclaim{Lemma C.1} For any $A\in \A_{l}(X_{r, R})$ with 
$l \geq 1$, there exists
a unique gauge $g\in \G_{l+1}(X_{r, R})$ with $g(x_0)=1$  such that
$\ti A=g^*(A)$ satisfies
$$\eqalign{d^*\ti A=0,\cr
A(\nu) =0 \hbox{ on }
\pa X_{r, R},\cr}\tag C.1
$$
where $\nu$ denotes the unit outer normal 
of $\pa X_{r, R}$. Moreover, we have 
$$\|A\|_{1, 2}\leq C\|F_A\|_{0, 2}$$
for a positive constant $C$ depending only on $Y, r $ and $R$.
\endproclaim
\demo{Proof}
The associated gauge fixing equation is 
$$\eqalign{d^*(g^{-1}dg)+d^*A & =0,\cr
g^{-1}dg(\nu)+A(\nu)& =0 \hbox{ on }
\pa X_{r, R}.\cr}\tag C.2 
$$
If we choose $g= e^f$, then the equation reduces to 
$$
\eqalign{d^*df+d^*A & =0,\cr
\frac{\pa f}{\pa \nu}+A(\nu)& =0 \hbox{ on }
\pa X_{r, R}.\cr}\tag C.3 
$$
It is clear that a solution   $f \in 
\O^0_{l+1}(X_{r, R})$ with $f(x_0) = 0$ exists.       

Now  assume that there are 
$g$  and $g_1$ satisfying (C.2) with 
$g(x_0)=g_1(x_0)=1$. Without loss of generality, we may assume that
$g_1=id$.  Taking  two copies of $X_{r, R}$ and gluing them along their 
common boundary, we obtain  the  Riemannian manifold 
$Y\times S^1$. The two copies  of $g$ then yield a solution $
g_0: Y \times S^1 \to S^1$ of the harmonic equation $d^*(g^{-1}dg)=0$. 
Clearly, for each $y\in Y$, $g_0(y,\cdot)$ is a map from $S^1$ to $S^1$
 with
degree zero, {\it i.e.} $\int_{\{y\}\times S^1}g_0^{-1}dg_0=0$. On the other 
hand, $g_0^{-1}dg_0$ defines an element $\o$ 
in $H^1_{de}$. Since $Y$ is a rational homology sphere, we deduce 
$\o=0$ and that $g_0$
is constant. Since $g(x_0)=0$, we infer  
$g\equiv 1$. 

 Next consider an $A$ satisfying (i) and (ii)
in Lemma C.1. As in the above argument, 
$A$ leads to a one form $A_0$ on $Y\times S^1$ 
satisfying $d^*A_0=0$. Since $Y$ is a rational homology sphere,
it is easy to see 
from the construction of $A_0$ that the harmonic part of $A_0$ in its 
Hodge decomposition zeor. Hence we have
$A_0=d^*\o$ for some two form $\o \in \O^2_{2}(Y\times S^1)$. 
Consequently, $d d^*\o= F_{A_0}$.
We deduce
$$\|A_0\|_{1, 2}= \|d^*\o\|_{1, 2}\le C\|F_{A_0}\|_{0, 2}.$$
This implies the desired estimate.
\qed
\enddemo

Next we prove transversality at reducible trajectories. 

\proclaim{Lemma C.2} Let $(\pi, \lmd)$ be a pair 
of good parameters such that $\n^2 H$ is sufficiently small. 
Choose $\d_-, d_+$ small enough (but positive)
in the set-up of Definition 4.4. Then the operator 
$\F_{p, q}$ for $p, q \in O$ at  a Seiberg-Witten 
trajectory is onto.
\endproclaim
\demo{Proof} Consider $\F= \F_{p, q}$ at a trajectory 
$(A_0, \P_0)$. By the proof of Lemma 4.12, we 
can assume $p=q=(a_0, 0)$ and $(A_0, \P_0)
\equiv p$. 
The formal 
adjoint $\F^*$ of $\F=\F_{p, q}$ with respect to the 
product (4.4) is given by 
$$\F^*(v)=-\dt v-\pmatrix
*db+df-\n^2 H(a_0)\cdot b \cr -\np_{a_0}\psi-\lmd\psi \cr
d^*b-\d_F' f\cr
\endpmatrix-{\d'(t)}v \tag C.4
$$
for $v=(\psi,b,f)$. The surjectivity of $\F$ 
is equivalent to the vanishing of the kernel of 
$\F^*$.  Let $v$ satisfy $\F^* v=0$. 
Then we have 
$$\dt\pmatrix b\cr f\cr\endpmatrix+
\pmatrix *db+df\cr d^*b\cr\endpmatrix+
\pmatrix {\d_F'}b+ \n^2H(a_0) \cdot b\cr 0\endpmatrix
=0. \tag C.5$$
We define the operator $L$ by 
$$L\pmatrix b\cr f\cr\endpmatrix=\pmatrix *db+df\cr d^*b\cr\endpmatrix.$$
$L$ is formally self-adjoint and  satisfies 
$L^2 = \D$, where $\D$ denotes the Hodge Laplacian. 
Let
$\{\xi_i=\pmatrix b_i\cr f_i\cr\endpmatrix\}$ be a complete $L^2$ 
orthonormal system of eignvectors of $L$ with 
$L\xi_i=\lmd_i\xi_i$. From the above discussion we deduce 
$$\D b_i=\lmd_i^2b_i,  \D f_i=\lmd_i^2 f_i.$$
Now we write 
$$\pmatrix b\cr f\cr\endpmatrix=\sum _{-\infty}^{+\infty}l_i(t)\xi_i.$$
Then it follows from C.5 that 
$$\sum l_i'(t)\xi_i+\lmd_il_i(t)\xi_i+\pmatrix 
{\d_F'}b + \n^2H(a_0) \cdot b \cr 0\endpmatrix=0.\tag C.6$$
Assume $\lmd_{j}=0$ for a $j$. Then there holds
$b_{j} =0$, for $Y$ is a rational homology sphere 
and hence supports no nontrivial harmonic 
1-form.  Consequently, $f_j$ is a nonzero 
constant. Then we deduce 
$l_{j}'(t) \equiv 0$, hence $l_j$ is a constant. But $(b, f)$ 
is $L^2$ integrable, which forces $l_j$ to be zero. We conclude 
that the above expansion of $(b ,f)$ does not contain terms 
with zero eigenvalue. Using the elementary 
arguments in e.g. $\Bz$ it is then easy to show that $(b, f)$ 
must vanish, provided that $\n^2H$, $\d_+$ and $\d_-$ have been chosen 
small enough.  Using the same arguments one also 
infers that $\psi$ vanishes. Thus $v=0$.  
\qed
\enddemo

\bigskip

{\it Aknowledgement: The first named author was supported 
by SFB 237 and the Leibniz program of DFG. He is grateful 
to Prof. J. Jost for providing these supports.}
\enddocument